# SEMIPARAMETRIC EFFICIENCY IN GMM MODELS WITH AUXILIARY DATA

By Xiaohong Chen,[1] Han Hong[2] and Alessandro Tarozzi

*New York University, Yale University, Stanford University and Duke University*

We study semiparametric efficiency bounds and efficient estimation of parameters defined through general moment restrictions with missing data. Identification relies on auxiliary data containing information about the distribution of the missing variables conditional on proxy variables that are observed in both the primary and the auxiliary database, when such distribution is common to the two data sets. The auxiliary sample can be independent of the primary sample, or can be a subset of it. For both cases, we derive bounds when the probability of missing data given the proxy variables is unknown, or known, or belongs to a correctly specified parametric family. We find that the conditional probability is not ancillary when the two samples are independent. For all cases, we discuss efficient semiparametric estimators. An estimator based on a conditional expectation projection is shown to require milder regularity conditions than one based on inverse probability weighting.

**1. Introduction.** Many empirical studies are complicated by the presence of missing data. One solution to the identification problem is based on the assumption that information on the true value of the variables in the data set of interest (the *primary* data set) can be recovered using *auxiliary* data sources under a conditional independence assumption. The key element of this identification strategy is that the distribution of the variables of interest is assumed to be independent of whether they belong to the primary or the auxiliary sample, conditional on a set of proxy variables, which are observed in both samples.

Received October 2005; revised May 2007.
[1]Supported in part by NSF Grants SES-03-18091 and SES-06-31613.
[2]Supported in part by NSF Grant SES-04-52143.

The first goal of this paper is to study semiparametric efficiency bounds of parameters defined through general nonlinear and over-identified moment conditions for missing data models under a conditional independence assumption. We provide semiparametric efficiency bounds for the cases when the propensity score is unknown, or is known, or belongs to a correctly specified parametric family. In our context, the propensity score is defined as the probability that one observation belongs to the subsample where only the proxy variables are observed. The auxiliary sample can be either a subset of the primary sample ("verify-in-sample" case) or independent of the primary sample ("verify-out-of-sample"). The former case is a special case of the MAR or CAR missing data structure where the missing variables are common to all subjects. Semiparametric efficiency bounds for this case are closely related to the results in Robins, Rotnitzky and Zhao (1994) when there is a single hierarchy in the case of monotone missing data patterns for a fixed set of instrument functions. [See also Robins and Rotnitzky (1995) and Chen and Breslow (2004).] We provide new results on semiparametric efficiency bounds for the "verify-out-of-sample" case. We find that while more information on the propensity score will not affect the asymptotic efficiency bounds for parameters defined in the verify-in-sample case [as shown in, e.g., Robins, Rotnitzky and Zhao (1994), Chen and Breslow (2004) and Hahn (1998)], it will improve the asymptotic efficiency for parameters defined in the verify-out-of-sample case. Our new efficiency bound results for the case when the parametric propensity is correctly specified should be useful in applied work because such an assumption is frequently adopted by empirical researchers.

The second goal of the paper is to develop two classes of sieve-based, Generalized Method of Moments (GMM) estimators that achieve the efficiency bounds for parameters defined under either the "verify-out-of-sample" or the "verify-in-sample" framework. Each estimator relies only on one nonparametric estimate; a conditional expectation projection based GMM (hereafter CEP-GMM) estimator only requires the nonparametric estimation of a conditional expectation, while an inverse probability weighting based GMM (hereafter IPW-GMM) estimator only needs a nonparametric estimate of the propensity score. We establish asymptotic normality and efficiency properties of both estimators under weaker regularity conditions than the existing ones in the literature. In particular, we allow for nonlinear and nonsmooth moment restrictions and for unbounded support of conditioning (or proxy) variables. The CEP-GMM estimator presents some advantages over the IPW-GMM estimator. First, its root-n asymptotic normality and efficiency can be derived without the strong assumption that the unknown propensity score is uniformly bounded away from zero and one. Second, the CEP-GMM estimator is characterized by a simple common format that achieves the relevant efficiency bound for all the cases we consider, regardless of whether the propensity score is unknown, or known or parametrically specified. Instead, the IPW-GMM estimator will be generally inefficient when the propensity score is known, or is parametrically estimated using a *correctly* specified parametric model; in such instances, different combinations of nonparametric and parametric estimates of the propensity score have to be specifically derived to achieve the efficiency bounds.

Our results can also be applied to the estimation of parametric nonlinear models with nonclassical measurement errors with validation data, a topic that has been studied in Carroll and Wand (1991), Sepanski and Carroll (1993), Carroll, Ruppert and Stefanski (1995), Lee and Sepanski (1995), Chen, Hong and Tamer (2005) among others.

Section 2 describes the model and presents the semiparametric efficiency bounds. Semiparametrically efficient CEP-GMM and IPW-GMM estimators are developed in Sections 3 and 4, respectively. In Section 5 we illustrate empirically the performance of the different estimators in the estimation of the distribution of private consumption in rural India in the presence of missing data. Section 6 concludes. All proofs are given in the Appendices.

**2. Semiparametric efficiency bounds.** Let $(X_i, Y_i, D_i)_{i=1}^n$ be an i.i.d. sample from $(X, Y, D)$, and denote $Z_i = (Y_i, X_i)$ where $Y_i$ is only observed when $D_i = 0$. We are interested in the estimation of parameters $\beta \in B$, a compact subset of $\mathcal{R}^{d_\beta}$, defined implicitly in terms of general nonlinear moment conditions. In the first (verify-out-of-sample) case such conditions are described by

$$(1) \qquad E[m(Z;\beta) \mid D = 1] = 0 \quad \text{if and only if} \quad \beta = \beta_0,$$

while in the second (verify-in-sample) case the condition is

$$(2) \qquad E[m(Z;\beta)] = 0 \quad \text{if and only if} \quad \beta = \beta_0,$$

where $Z = (Y, X)$ and $m(\cdot;\beta)$ is a set of functions with dimension $d_m \geq d_\beta$.

In other words, under case (1) $Y$ is always missing in the *primary* data set $(D = 1)$, which is a random sample from the population of interest, while an independent *auxiliary* sample (where $D = 0$) will serve the purpose of ensuring the identification of parameters that would not be identified by the primary data set alone. Under case (2), the auxiliary sample is instead a subset of the entire primary sample.

In this section we present the semiparametric efficiency bound for the estimation of $\beta$ implicitly defined by either moment conditions (1) or (2). In this paper $\beta$ is typically used to denote an arbitrary value in the parameter space, but to save notation in this section $\beta$ is also used as the true parameter value $\beta_0$. Define

$$(3) \qquad \mathcal{E}(X;\beta) = E[m(Z;\beta) \mid X]$$

to be the conditional expectation of the moment conditions given $X$, and define

(4) $\quad V(m(Z;\beta)|X) = E[m(Z;\beta)m(Z;\beta)' \mid X] - \mathcal{E}(X;\beta)\mathcal{E}(X;\beta)'$

to be the conditional variance of the moment conditions given $X$. In addition, define

$$p = \Pr(D=1) \quad \text{and} \quad p(X) = \Pr(D=1 \mid X),$$

$$\mathcal{J}_\beta^1 = \frac{\partial}{\partial \beta} E[m(Z;\beta) \mid D=1] \quad \text{and} \quad \mathcal{J}_\beta^2 = \frac{\partial}{\partial \beta} E[m(Z;\beta)].$$

ASSUMPTION 1. (i) Both $\mathcal{J}_\beta^1$ and $\mathcal{J}_\beta^2$ have full column rank equal to $d_\beta$; (ii) The data $(X_i, Y_i, D_i)_{i=1}^n$ is an i.i.d. sample from $(X, Y, D)$; (iii) $p \in (0,1)$, $p(X) \in (0,1)$.

Notice that in both cases (1) and (2) the moment conditions are assumed to hold in the *primary* sample in which some information is missing. Identification is possible because of the access to an auxiliary data set ($D=0$) which contains both $Y$ and a set of proxy variables $X$ that are also potentially of interest, if the following fundamental conditional independence assumption holds:

ASSUMPTION 2. $Y \perp D \mid X$.

Conditional independence assumptions have been used extensively in econometrics and statistics to achieve identification with missing data. Examples include inference in models with attrition or nonresponse [e.g., Little and Rubin (2002), Robins and Rotnitzky (1995), Rotnitzky and Robins (1995), Wooldridge (2002), Wooldridge (2003)], the estimation of treatment effects [see e.g., the references surveyed in Heckman, LaLonde and Smith (1999)], the recovery of comparability over time of statistics calculated using data collected with different methodology [e.g., Clogg, Rubin, Schenker, Schultz and Weidman (1991), Schenker (2003), Tarozzi (2007)].

Under case (2), Assumption 1 would be satisfied if, for instance, the probability of validating a given observation only depends on $X$. In case (1), Assumption 2 requires that the sampling scheme used to create the auxiliary sample depends only on $X$. If a simple random subset of the primary data is validated, $p(X)$ is a constant and the auxiliary data set is characterized by the same distribution of $(Y, X)$ as the primary data set, and Assumption 2 is easily seen satisfied. In this case, which is common in the statistics literature, the auxiliary data set is usually called a *validation* data set. A stratified sample satisfying Assumption 2 in model

(2) can also be produced through a two-stage sampling design using a finite number of strata [see e.g., Breslow, Robins and Wellner (2000) and Breslow, McNeney and Wellner (2003)], in which case the only variable $X$ that is observed for all sampled observations is a discrete stratum indicator. In the following, the "regular estimators" are defined according to Begun, Hall, Huang and Wellner (1983) and Ibragimov and Has'minskii (1981).

THEOREM 1. *Let $\beta$ be defined by the moment conditions (1) or (2). Under Assumptions 1–2, the asymptotic variance lower bound for all regular estimators of $\beta$ is*

$$(\mathcal{J}_\beta' \Omega_\beta^{-1} \mathcal{J}_\beta)^{-1} \qquad \text{for some } \mathcal{J}_\beta \text{ and some positive definite } \Omega_\beta,$$

*where, for the moment condition (1), $\mathcal{J}_\beta = \mathcal{J}_\beta^1$ and $\Omega_\beta = \Omega_\beta^1$:*

$$\Omega_\beta^1 = E\left[\frac{p(X)^2}{p^2(1-p(X))}V(m(Z;\beta) \mid X) + \frac{p(X)}{p^2}\mathcal{E}(X;\beta)\mathcal{E}(X;\beta)'\right];$$

*and for the moment condition (2), $\mathcal{J}_\beta = \mathcal{J}_\beta^2$ and $\Omega_\beta = \Omega_\beta^2$:*

$$\Omega_\beta^2 = E\left[\frac{1}{1-p(X)}V[m(Z;\beta) \mid X] + \mathcal{E}(X;\beta)\mathcal{E}(X;\beta)'\right].$$

In Appendix A we present explicit expressions for the efficient score functions corresponding to the asymptotic variance lower bounds in Theorem 1 as well as in the following Theorems 2 and 3.

2.1. *Information content of the propensity score.* It is interesting to analyze whether the knowledge of the propensity score $p(X)$ decreases the semiparametric efficiency bounds for the parameters $\beta$.

THEOREM 2. *Let $\beta$ be defined by the moment conditions (1) or (2). Under Assumptions 1–2, if $p(X)$ is known, then the asymptotic variance lower bound for all regular estimators of $\beta$ is*

$$(\mathcal{J}_\beta' \tilde{\Omega}_\beta^{-1} \mathcal{J}_\beta)^{-1} \qquad \text{for some } \mathcal{J}_\beta \text{ and some positive definite } \tilde{\Omega}_\beta,$$

*where, for the moment condition (1), $\mathcal{J}_\beta = \mathcal{J}_\beta^1$ and $\tilde{\Omega}_\beta = \tilde{\Omega}_\beta^1$:*

$$\tilde{\Omega}_\beta^1 = E\left[\frac{p(X)^2}{p^2(1-p(X))}V(m(Z;\beta) \mid X) + \frac{p(X)^2}{p^2}\mathcal{E}(X;\beta)\mathcal{E}(X;\beta)'\right];$$

*and for the moment condition (2), $\mathcal{J}_\beta = \mathcal{J}_\beta^2$ and $\tilde{\Omega}_\beta = \Omega_\beta^2$ given in Theorem 1.*

In other words, knowledge of $p(X)$ reduces the semiparametric efficient variance bound for $\beta$ under the "verify-out-of-sample" case, but it does not under the "verify-in-sample" case. The following argument provides an intuition for this result. When (2) holds, $\beta$ is defined through the relation

$$\iint m(y,x;\beta)f(y \mid x)\,dy f(x)\,dx = 0.$$

The propensity score $p(X)$ does not enter the definition of $\beta$, therefore its knowledge should not affect the variance bound for $\beta$. However, the relation that identifies $\beta$ when (1) holds clearly depends on $p(X)$:

$$\iint m(y,x;\beta)p(x)f(y \mid x)\,dy f(x)\,dx = 0.$$

REMARK 1. A special case of Theorem 2 is when $p(X)$ is a constant $p$. In this case, the auxiliary sample is also called a validation sample and is drawn randomly from the same population as the primary sample, so that $Y, X \perp D$ [Carroll and Wand (1991), Sepanski and Carroll (1993), Lee and Sepanski (1995)]. In such case it is then easy to see that the two efficiency bounds given in Theorem 2 become identical.

Another interesting question is what is the efficiency bound for the estimation of $\beta$ defined by moment condition (1) if the propensity score is unknown but is assumed to belong to a correctly specified parametric family, so that $p(X) = p(X;\gamma)$. Let $p_\gamma(X) = \partial p(X;\gamma)/\partial \gamma$, and define the score function for $\gamma$ as $S_\gamma = S_\gamma(D,X) = \frac{D - p(X;\gamma)}{p(X;\gamma)(1-p(X;\gamma))} p_\gamma(X)$.

THEOREM 3. *Let $\beta$ be defined by the moment conditions* (1). *Under Assumptions* 1–2, *if $p(X) = p(X;\gamma)$ and $E[S_\gamma(D,X)S_\gamma(D,X)']$ is positive definite, then the asymptotic variance lower bound for all regular estimators of $\beta$ is $(\mathcal{J}_\beta' \tilde{\Omega}_\beta^{-1} \mathcal{J}_\beta)^{-1}$, where $\mathcal{J}_\beta = \mathcal{J}_\beta^1$, $\tilde{\Omega}_\beta^1$ is given in Theorem* 2 *and*

$$\tilde{\Omega}_\beta = \tilde{\Omega}_\beta^1 + \left[E\frac{\mathcal{E}(X;\beta)p_\gamma(X)'}{p}\right][ES_\gamma S_\gamma']^{-1}\left[E\frac{p_\gamma(X)\mathcal{E}(X;\beta)'}{p}\right].$$

This variance bound is clearly larger than $\tilde{\Omega}_\beta^1$ stated in Theorem 2, but it is smaller than the bound in Theorem 1. This latter result can be verified noting first that the bound in Theorem 3 corresponds to the variance of the following influence function:

$$\frac{(1-D)p(X)}{p(1-p(X))}(m(Z;\beta) - \mathcal{E}(X;\beta))$$
$$+ \text{Proj}\left(\frac{\mathcal{E}(X;\beta)}{p}(D - p(X))\Big| S_\gamma(D,X)\right) + \frac{p(X)\mathcal{E}(X;\beta)}{p},$$

where we use $\text{Proj}(Z_1|Z_2)$ to denote the population least squares projection of a random variable $Z_1$ onto the linear space spanned by $Z_2$. The conclusion follows noting that the variance bound stated in Theorem 1 for moment condition (1) is instead the variance of the following influence function

$$\frac{1}{p}D\mathcal{E}(X;\beta) + \frac{(1-D)p(X)}{p(1-p(X))}[m(Z;\beta) - \mathcal{E}(X;\beta)],$$

whose corresponding variance is larger.

Our results for GMM models complement and extend the finding in the program evaluation literature that knowing the propensity score decreases the efficient variance bound for the estimation of the average effect of treatment on the treated, while the propensity score is ancillary for the average treatment effect parameter [Hahn (1998)].

**3. CEP-GMM estimation.** In this section, we consider a first class of semiparametrically efficient estimators based on a conditional expectation projection (CEP) method. If Assumption 2 holds, identification follows by noting that, under case (1)

$$E[m(Z;\beta) \mid D = 1] = \int E[m(Z;\beta) \mid x, D = 0]f(x \mid D = 1)\,dx,$$

while under case (2),

$$E[m(Z;\beta)] = \int E[m(Z;\beta) \mid x, D = 0]f(x)\,dx.$$

Therefore, estimation of the parameters of interest can proceed by first estimating $E[m(Z;\beta)|x, D = 0]$ nonparametrically from the auxiliary sample, and then integrating the conditional expectation against the distribution of $x$ in the primary sample.

3.1. *Efficient estimation with unknown propensity score.* In the following, we use subscripts $p$ and $a$ to refer respectively to observations belonging to the primary and to the auxiliary sample. Let $n_p$ be the size of the primary sample and $n_a$ be the size of the auxiliary sample. Observations in the primary sample are indexed by $i = 1, \ldots, n_p$. Observations in the auxiliary sample are indexed by $j = 1, \ldots, n_a$. Under moment condition (1) (verify-out-of-sample case), $n = n_p + n_a$. Under moment condition (2) (verify-in-sample case), $n = n_p$. Let $\hat{\mathcal{E}}(X;\beta)$ denote a nonparametric estimate of $\mathcal{E}(X;\beta)$ using the auxiliary sample. Chen, Hong and Tamer (2005) (hereafter CHT) used a sieve Least Squares (LS) estimator. Let $\{q_l(X), l = 1, 2, \ldots\}$ denote a sequence of known basis functions that can approximate any square-measurable function of $X$ arbitrarily well. Also let

$$q^{k(n_a)}(X) = (q_1(X), \ldots, q_{k(n_a)}(X))'$$

and

$$Q_a = (q^{k(n_a)}(X_{a1}), \ldots, q^{k(n_a)}(X_{an_a}))'$$

for some integer $k(n_a)$, with $k(n_a) \to \infty$ and $k(n_a)/n \to 0$ when $n \to \infty$. Then for each given $\beta$, the sieve LS estimator of $\mathcal{E}(X;\beta)$ is

$$\hat{\mathcal{E}}(X;\beta) = \sum_{j=1}^{n_a} m(Z_{aj};\beta) q^{k(n_a)}(X_{aj})(Q_a'Q_a)^{-1} q^{k(n_a)}(X).$$

A generalized method of moment estimator for $\beta_0$ can then be defined as

$$(5) \qquad \hat{\beta} = \arg\min_{\beta \in B} \left( \frac{1}{n_p} \sum_{i=1}^{n_p} \hat{\mathcal{E}}(X_{pi};\beta) \right)' \hat{W} \left( \frac{1}{n_p} \sum_{i=1}^{n_p} \hat{\mathcal{E}}(X_{pi};\beta) \right).$$

The $\sqrt{n}$-consistency and asymptotic normality of this CEP-GMM estimator have been established in CHT. Following the proof of their claim (A.2), we have the following asymptotic representation:

$$\frac{\sqrt{n}}{n_p} \sum_{i=1}^{n_p} \hat{\mathcal{E}}(X_{pi};\beta_0) = \frac{\sqrt{n}}{n_p} \sum_{i=1}^{n_p} \mathcal{E}(X_{pi};\beta_0)$$

$$+ \frac{\sqrt{n}}{n_a} \sum_{j=1}^{n_a} \frac{f_{X_p}(X_{aj})}{f(X_{aj}|D=0)} [m(Z_{aj};\beta_0) - \mathcal{E}(X_{aj};\beta_0)]$$

$$+ o_p(1),$$

where we use $f_{X_p}(X)$ to denote the density of $X$ in the primary data set, and $o_p(1)$ represents a term that converges to 0 in probability.

When moment condition (1) holds, $n = n_p + n_a$, $f_{X_p}(X) = f(X \mid D = 1)$ and

$$\frac{f_{X_p}(X)}{f(X \mid D = 0)} = \frac{(1-p)p(X)}{p(1-p(X))}.$$

In this case we can also write the influence function for $\frac{\sqrt{n}}{n_p} \sum_{i=1}^{n_p} \hat{\mathcal{E}}(X_{pi};\beta_0)$ as

$$\frac{1}{\sqrt{n}} \sum_{i=1}^{n} \left\{ \frac{1}{p} D_i \mathcal{E}(X_i;\beta_0) \right.$$

$$(6) \qquad \qquad \left. + (1-D_i) \frac{p(X_i)}{p(1-p(X_i))} [m(Z_i;\beta_0) - \mathcal{E}(X_i;\beta_0)] \right\}$$

$$+ o_p(1).$$

The proof of Theorem 1 shows that the two terms in the influence function correspond to the two components of the efficient influence function that contain information about $f(X|D=1)$ and $f(Y \mid X)$, respectively. These two terms are orthogonal to each other, so that

$$Avar\left(\frac{\sqrt{n}}{n_p}\sum_{i=1}^{n_p}\hat{\mathcal{E}}(X_{pi};\beta_0)\right) = \Omega_\beta^1,$$

where $\Omega_\beta^1$ is given in Theorem 1.

When moment condition (2) holds, $f_{X_p}(X) = f(X)$, $n_p = n$ and

$$\frac{f_{X_p}(X)}{f(X \mid D=0)} = \frac{1-p}{1-p(X)}.$$

The influence function for $\frac{\sqrt{n}}{n_p}\sum_{i=1}^{n_p}\hat{\mathcal{E}}(X_{pi};\beta_0)$ can then be written as

$$\frac{1}{\sqrt{n}}\sum_{i=1}^{n}\left\{\mathcal{E}(X_i;\beta_0) + (1-D_i)\frac{1}{1-p(X_i)}[m(Z_i;\beta_0) - \mathcal{E}(X_i;\beta_0)]\right\}$$
(7)
$$+ o_p(1).$$

The two terms in the influence function correspond to the two components of the projected efficiency influence function that contain information about $f(X)$ and $f(Y \mid X)$ in the proof of Theorem 1. The orthogonality between these two terms implies that

$$Avar\left(\frac{\sqrt{n}}{n_p}\sum_{i=1}^{n_p}\hat{\mathcal{E}}(X_{pi};\beta_0)\right) = \Omega_\beta^2,$$

where $\Omega_\beta^2$ is given in Theorem 1. The semiparametric efficiency bounds given in Theorem 1 are then achieved by an optimally weighted GMM estimator $\hat{\beta}$ for $\beta_0$ that uses a weighting matrix $\hat{W} = \Omega_\beta^{-1} + o_p(1)$.

THEOREM 4. *Let $\hat{\beta}$ be the CEP-GMM estimator given in (5). Under Assumptions 1–2, and Assumptions 3–5 of CHT, we have $\sqrt{n}(\hat{\beta} - \beta_0) \Rightarrow \mathcal{N}(0,V)$, with $V = (\mathcal{J}_\beta' W \mathcal{J}_\beta)^{-1}\mathcal{J}_\beta' W \Omega_\beta W \mathcal{J}_\beta(\mathcal{J}_\beta' W \mathcal{J}_\beta)^{-1}$, where $\Omega_\beta$ is given in Theorem 1. Furthermore, if $W = \Omega_\beta^{-1}$, then $\sqrt{n}(\hat{\beta} - \beta_0) \Rightarrow \mathcal{N}(0,V_0)$, with $V_0 = (\mathcal{J}_\beta' \Omega_\beta^{-1} \mathcal{J}_\beta)^{-1}$, where $\mathcal{J}_\beta = \mathcal{J}_\beta^1$ and $\Omega_\beta = \Omega_\beta^1$ under moment condition (1), and $\mathcal{J}_\beta = \mathcal{J}_\beta^2$ and $\Omega_\beta = \Omega_\beta^2$ under moment condition (2).*

CHT derive the root-$n$ consistency and normality of the CEP-GMM estimator. Theorem 4 says that their estimator $\hat{\beta}$ is also semiparametrically efficient. The proof of Theorem 4 follows directly from that of Theorem 2 in CHT, who also provide simple consistent estimators of $V$ and $V_0$. In the working paper version of this article, we have stated Assumptions 3–5 of CHT in terms of the notations of this paper. These assumptions are mild regularity conditions. In particular, they allow for (i) nonsmooth $m(Z;\beta)$, such as quantile-based moment functions; (ii) the support of the conditioning (proxy) variable $X$ could be unbounded; (iii) the propensity score function $p(X)$ does not need to be uniformly bounded away from zero and one. Recall that in the program evaluation literature such as in Hirano, Imbens and Ridder (2003), the stronger condition $0 < \underline{p} \leq p(x) \leq \overline{p} < 1$ is typically imposed for root-$n$ asymptotically normal and efficient estimation of $\beta_0$.

3.2. *CEP estimation with parametric or known propensity score.* Suppose now that the propensity score $p(X)$ is correctly parameterized as $p(X;\gamma)$ up to a finite-dimensional unknown parameter $\gamma$. Theorems 2 and 4 show that the optimally weighted CEP-GMM estimator defined in (5) still achieves the semiparametric efficiency bound for $\beta$ defined by moment condition (2). However, according to Theorems 3 and 4, such an estimator is no longer efficient for $\beta$ defined through moment condition (1).

Rewriting moment condition (1) as $E[\mathcal{E}(X_i;\beta_0)\frac{p(X_i)}{p}] = 0$, we can again construct an efficient estimator for $\beta_0$ based on the sieve estimate $\hat{\mathcal{E}}(X;\beta)$ and the correctly specified parametric form $p(X;\gamma)$. In particular, the optimally weighted GMM estimator using the following sample moment condition will achieve the efficiency bound in Theorem 3 for $\beta$ defined through (1):

$$(8) \qquad \frac{1}{n}\sum_{i=1}^{n}\hat{\mathcal{E}}(X_i;\beta)\frac{p(X_i;\hat{\gamma})}{\hat{p}},$$

where $\hat{p} = \frac{n_p}{n}$ and $\hat{\gamma}$ is the parametric MLE that solves the score equation for $\gamma$:

$$\frac{1}{n}\sum_{i=1}^{n}S_{\hat{\gamma}}(D_i,X_i) = \frac{1}{n}\sum_{i=1}^{n}\frac{D_i - p(X_i;\hat{\gamma})}{p(X_i;\hat{\gamma})(1-p(X_i;\hat{\gamma}))}p_{\hat{\gamma}}(X_i) = 0.$$

THEOREM 5. *Let $p(X;\gamma)$ be the parametric propensity score function known up to the parameters $\gamma$ and let $E[S_{\gamma_0}(D,X)S_{\gamma_0}(D,X)']$ be positive definite. Let $\beta_0$ satisfy the moment condition* (1) *and $\widehat{\beta}$ be its CEP-GMM estimator using the sample moment* (8). *Under Assumptions* 1–2 *and Assumptions* 3–5 *of CHT, we have $\sqrt{n}(\widehat{\beta} - \beta_0) \Rightarrow \mathcal{N}(0,V)$, with*

$$V = (\mathcal{J}_{\beta}^{1\prime}W\mathcal{J}_{\beta}^{1})^{-1}\mathcal{J}_{\beta}^{1\prime}W\tilde{\Omega}_{\beta}W\mathcal{J}_{\beta}^{1}(\mathcal{J}_{\beta}^{1\prime}W\mathcal{J}_{\beta}^{1})^{-1},$$

*where $\tilde{\Omega}_{\beta}$ is given in Theorem* 3. *Further, if $W = \tilde{\Omega}_{\beta}^{-1}$, then $\sqrt{n}(\widehat{\beta} - \beta_0) \Rightarrow \mathcal{N}(0,V_0)$, where $V_0 = (\mathcal{J}_{\beta}^{1\prime}\tilde{\Omega}_{\beta}^{-1}\mathcal{J}_{\beta}^{1})^{-1}$ is the efficiency variance bound given in Theorem* 3.

The proof of this theorem is similar to that of Theorem 2 in CHT and is thus omitted. The influence function representation of (8) is stated in the working paper version of this article.

We remark that even when a parametric assumption is being made about the propensity score $p(X; \gamma)$ [in fact even if in addition $f(Y)$ is assumed to be a parametric likelihood], the inference about $\beta$ is still semiparametric. This is because the marginal density $f(X)$ is still nonparametric and contains semiparametric information about $\beta$. This explains why nonparametric estimation is still needed to achieve the bound for $\beta$.

The case where the propensity score is fully known can be considered a special case of parametric propensity score where the parameters are known. In this case, the efficient moment condition is as in (8) after replacing $p(X_j; \hat{\gamma})$ with the known $p(X_j)$.

REMARK 2. When the auxiliary data set is a validation data set, for example, $p(X) = p$, the parameters $\beta$ defined by both moment conditions (1) and (2) coincide. Therefore, the CEP-GMM estimator defined in (5) when we take $n_p = n$ and the summation to be over the all observations will achieve semiparametric efficiency.

**4. IPW-GMM estimation.** An alternative estimation method for $\beta$ is the inverse probability weighting based GMM (IPW-GMM). Several authors have considered inverse probability weighting paired with a conditional independence assumption for estimation in presence of missing information. Recent examples include parametric IPW as in Robins, Mark and Newey (1992), Wooldridge (2002), Wooldridge (2003) and Tarozzi (2007), for missing data models, and nonparametric inverse probability weighting as in Hirano, Imbens and Ridder (2003) for the case of mean treatment effect analysis. In this section, we extend existing results and first show that the optimally weighted IPW-GMM estimator of $\beta$ is semiparametrically efficient when the propensity score is unknown. The same estimator, however, will be generally inefficient when the propensity score is known or belongs to a correctly specified parametric family; combinations of nonparametric and known or parametric estimated propensity scores are needed to achieve the semiparametric efficiency bounds for these cases.

4.1. *Efficient estimation with unknown propensity score.* The IPW-GMM method uses the fact that under Assumption 2, moment condition (1) can be rewritten as

(9) $$E[m(Z;\beta) \mid D = 1] = E\left[m(Z;\beta)\frac{p(X)(1-p)}{(1-p(X))p} \,\Big|\, D = 0\right],$$

while moment condition (2) is equivalent to

$$E[m(Z;\beta)] = E\left[m(Z;\beta)\frac{1-p}{1-p(X)} \,\Big|\, D=0\right]. \tag{10}$$

Let $\hat{p}(X)$ be a consistent estimate of the true propensity score. Then we can estimate $\beta_0$ defined by case (1) using GMM with the following sample moment:

$$\sqrt{n}\frac{1}{n_a}\sum_{j=1}^{n_a} m(Z_j;\beta)\frac{\hat{p}(X_j)}{1-\hat{p}(X_j)}\frac{1-\hat{p}}{\hat{p}}, \tag{11}$$

and estimate $\beta_0$ defined by case (2) using GMM with the following sample moment:

$$\sqrt{n}\frac{1}{n_a}\sum_{j=1}^{n_a} m(Z_j;\beta)\frac{1-\hat{p}}{1-\hat{p}(X_j)}. \tag{12}$$

The inverse probability weighting approach is considered semiparametric when $\hat{p}(X)$ is estimated nonparametrically. In this case, it can be shown that the sample moment (11) evaluated at $\beta_0$ is asymptotically equivalent to

$$\frac{1}{p}\frac{1}{\sqrt{n}}\sum_{i=1}^{n}\left[(1-D_i)m(Z_i;\beta_0)\frac{p(X_i)}{1-p(X_i)} + \mathcal{E}(X_i;\beta_0)\frac{D_i-p(X_i)}{1-p(X_i)}\right] \\ + o_p(1). \tag{13}$$

The two components of this influence function are negatively correlated. Because of this, the asymptotic variance might be smaller than that of the estimator of $\beta_0$ based on moment condition (11) with the known $p(X)$. Simple manipulations are sufficient to show that (13) is identical to the influence function in (6) whose corresponding asymptotic variance is $\Omega_\beta^1$ given in Theorem 1. An optimally weighted GMM estimator for $\beta_0$ defined by case (1) using sample moment (11) achieves then the semiparametric efficiency bound stated in Theorem 1.

The influence function representation for sample moment (12) can be calculated as

$$\frac{1}{\sqrt{n}}\sum_{i=1}^{n}\left[(1-D_i)m(Z_i;\beta_0)\frac{1}{1-p(X_i)} + \mathcal{E}(X_i;\beta_0)\frac{D_i-p(X_i)}{1-p(X_i)}\right] + o_p(1),$$

whose two components are again negatively correlated. However, it is again simple to show that this influence function representation is identical to the one in (7). Hence, an optimally weighted GMM estimator for $\beta_0$ defined by case (2) using sample moment (12) achieves the bound for case (1) stated in Theorem 1.

In this subsection, to emphasize that the true propensity score function is unknown and has to be estimated nonparametrically, we use $p_o(x) \equiv E(D \mid X = x)$ to indicate the true propensity score and $p(x)$ to denote any candidate function. [Note that to save notations in the rest of the main text $p(x)$ denotes the true propensity score.] Let $\widehat{p}(\cdot)$ be a sieve estimator of $p_o(x)$ that uses the combined sample $\{(D_i, X_i) : i = 1, \ldots, n\}$. Let $\{Z_{ai} = (Y_{ai}, X_{ai}) : i = 1, \ldots, n_a\}$ be the auxiliary (i.e., $D = 0$) data set. We define the IPW-GMM estimator $\widehat{\beta}$ for moment condition (1) as

$$(14) \quad \widehat{\beta} = \arg\min_{\beta \in B} \left( \frac{1}{n_a} \sum_{i=1}^{n_a} m(Z_{ai}; \beta) \frac{\widehat{p}(X_{ai})}{1 - \widehat{p}(X_{ai})} \right)' \widehat{W}$$
$$\times \left( \frac{1}{n_a} \sum_{i=1}^{n_a} m(Z_{ai}; \beta) \frac{\widehat{p}(X_{ai})}{1 - \widehat{p}(X_{ai})} \right)$$

and the IPW-GMM estimator $\widehat{\beta}$ for moment condition (2) as

$$(15) \quad \widehat{\beta} = \arg\min_{\beta \in B} \left( \frac{1}{n_a} \sum_{i=1}^{n_a} m(Z_{ai}; \beta) \frac{1}{1 - \widehat{p}(X_{ai})} \right)' \widehat{W}$$
$$\times \left( \frac{1}{n_a} \sum_{i=1}^{n_a} m(Z_{ai}; \beta) \frac{1}{1 - \widehat{p}(X_{ai})} \right).$$

There are two popular sieve nonparametric estimators of $p_o(\cdot)$:

(i) A sieve Least Squares (LS) estimator $\widehat{p}_{\text{ls}}(x)$ as in Hahn (1998):

$$\widehat{p}_{\text{ls}} = \arg\min_{p(\cdot) \in \mathcal{H}_n} \frac{1}{n} \sum_{i=1}^{n} (D_i - p(X_i))^2 / 2,$$

$$\mathcal{H}_n = \left\{ h(x) = q^{k_n}(x)' \pi = \sum_{j=1}^{k_n} q_j(x) \pi_j \right\} \quad \text{for some known basis } (q_j)_{j=1}^{\infty}.$$

In the Appendix we establish the consistency and convergence rate of $\widehat{p}_{\text{ls}}(x)$ under the assumption that the variables in $X$ have unbounded support.

(ii) A sieve Maximum Likelihood (ML) estimator $\widehat{p}_{\text{mle}}(x)$ as in Hirano, Imbens and Ridder (2003):

$$\widehat{p}_{\text{mle}} = \arg\max_{p(\cdot) \in \mathcal{H}_n} \frac{1}{n} \sum_{i=1}^{n} \{D_i \log[p(X_i)] + (1 - D_i) \log[1 - p(X_i)]\},$$

$$\mathcal{H}_n = \{h(x) = [A^{k_n}(x)' \pi]^2\} \text{ or } \{h(x) = \exp(A^{k_n}(x)' \pi)\}.$$

Before we present the large sample properties of the IPW-GMM estimator, we need to introduce some notations and assumptions. Let the support of $X$ be $\mathcal{X} = \mathcal{R}^{d_x}$. We could use more complicated notations and let $\mathcal{X} = \mathcal{X}_{\mathrm{c}} \times \mathcal{X}_{\mathrm{dc}}$, with $\mathcal{X}_{\mathrm{c}}$ being the support of the continuous variables and $\mathcal{X}_{\mathrm{dc}}$ the support of the finitely many discrete variables. Further we could decompose $\mathcal{X}_{\mathrm{c}} = \mathcal{X}_{\mathrm{c}1} \times \mathcal{X}_{\mathrm{c}2}$ with $\mathcal{X}_{\mathrm{c}1} = \mathcal{R}^{d_{x,1}}$ and $\mathcal{X}_{c2}$ being a compact and connected subset of $\mathcal{R}^{d_{x,2}}$. Then, under simple and usual modifications of the assumptions, the large sample results stated below would remain valid. To avoid tedious notation yet to allow for some unbounded support elements of $X$, we assume $\mathcal{X} = \mathcal{X}_{\mathrm{c}} = \mathcal{R}^{d_x}$. For any $1 \times d_x$ vector $\mathbf{a} = (a_1, \ldots, a_{d_x})$ of nonnegative integers, we write $|\mathbf{a}| = \sum_{k=1}^{d_x} a_k$, and for any $x = (x_1, \ldots, x_{d_x})' \in \mathcal{X}$, we denote the $|\mathbf{a}|$th derivative of a function $h : \mathcal{X} \to \mathcal{R}$ as

$$\nabla^{\mathbf{a}} h(x) = \frac{\partial^{|\mathbf{a}|}}{\partial x_1^{a_1} \cdots \partial x_{d_x}^{a_{d_x}}} h(x).$$

For some $\gamma > 0$, let $\underline{\gamma}$ be the largest integer smaller than $\gamma$, and let $\Lambda^\gamma(\mathcal{X})$ denote a Hölder space with smoothness $\gamma$, that is, a space of functions $h : \mathcal{X} \to \mathcal{R}$ which have up to $\underline{\gamma}$ continuous derivatives, and the highest ($\underline{\gamma}$th) derivatives are Hölder continuous with the Hölder exponent $\gamma - \underline{\gamma} \in (0, 1]$. The Hölder space becomes a Banach space when endowed with the Hölder norm:

$$\|h\|_{\Lambda^\gamma} = \sup_x |h(x)| + \max_{|\mathbf{a}|=\underline{\gamma}} \sup_{x \neq \overline{x}} \frac{|\nabla^{\mathbf{a}} h(x) - \nabla^{\mathbf{a}} h(\overline{x})|}{\sqrt{(x-\overline{x})'(x-\overline{x})}^{\gamma-\underline{\gamma}}} < \infty.$$

We call $\Lambda^\gamma_c(\mathcal{X}) \equiv \{h \in \Lambda^\gamma(\mathcal{X}) \| h \|_{\Lambda^\gamma} \leq c < \infty\}$ a Hölder ball (with radius $c$).

Define a weighted sup-norm $\|g\|_{\infty,\omega} \equiv \sup_{x \in \mathcal{X}} |g(x)[1+|x|^2]^{-\omega/2}|$ for some $\omega > 0$. Denote $\Pi_{\infty n} g$ as the projection of $g$ onto the sieve space $\mathcal{H}_n$ under the norm $\|\cdot\|_{\infty,\omega}$. Let $f_{X_a}(x) = f_{X|D=0}(x)$ and $E_a(\cdot) = E(\cdot|D=0)$.

ASSUMPTION 3. Let $\widehat{W} - W = o_p(1)$ for a positive semidefinite matrix $W$, and the followings hold: (1) $p_o(\cdot)$ belongs to a Hölder ball $\mathcal{H} = \{p(\cdot) \in \Lambda^\gamma_c(\mathcal{X}) : 0 < \underline{p} \leq p(x) \leq \overline{p} < 1\}$ for some $\gamma > 0$; (2) $\int (1+|x|^2)^\omega f_X(x)\, dx < \infty$ for some $\omega > 0$; (3) there is a function $b(\cdot)$ such that $b(\delta) \to 0$ as $\delta \to 0$ and $E_a[\sup_{\|\beta - \widetilde{\beta}\| < \delta} \|m(Z_i; \beta) - m(Z_i, \widetilde{\beta})\|^2] \leq b(\delta)$ for all small positive value $\delta$; (4) $E_a[\sup_{\beta \in B} \|m(Z_i; \beta)\|^2] < \infty$; (5) for any $h \in \mathcal{H}$, there is a sequence $\Pi_{\infty n} h \in \mathcal{H}_n$ such that $\|h - \Pi_{\infty n} h\|_{\infty,\omega} = o(1)$.

THEOREM 6. Let $\widehat{\beta}$ be the IPW-GMM estimator given in (14) or (15). Under Assumptions 1, 2 and 3, if $\frac{k_n}{n} \to 0$, $k_n \to \infty$, then $\widehat{\beta} - \beta_0 = o_p(1)$.

Let $E(\cdot) = \int (\cdot) f_X(x)\, dx$, $\|h\|_2 = \sqrt{\int h(x)^2 f_X(x)\, dx}$, and $\Pi_{2n} h$ be the projection of $h$ onto the closed linear span of $q^{k_n}(x) = (q_1(x), \ldots, q_{k_n}(x))'$ under the norm $\|\cdot\|_2$. We need the following additional assumptions to obtain asymptotic normality.

ASSUMPTION 4. Let $\beta_0 \in int(B)$, $E[\frac{p_o(X)}{1-p_o(X)}\mathcal{E}(X;\beta_0)\mathcal{E}(X;\beta_0)']$ be positive definite, and the followings hold: (1) Assumptions 3.1 and 3.2 are satisfied with $\gamma > d_x/2$ and $\omega > \gamma$; (2) There exist a constant $\epsilon \in (0,1]$ and a small $\delta_0 > 0$ such that

$$E_a\left[\sup_{\|\beta-\widetilde{\beta}\|<\delta} \|m(Z_i;\beta) - m(Z_i,\widetilde{\beta})\|^2\right] \leq const.\delta^\epsilon$$

for any small positive value $\delta \leq \delta_0$; (3) $E_a[\sup_{\beta \in B: \|\beta-\beta_0\|\leq\delta_0} \|m(Z_i;\beta)\|^2(1+|X_i|^2)^\omega] < \infty$ for some small $\delta_0 > 0$; (4) $E[\|\frac{\partial \mathcal{E}(X;\beta_0)}{\partial \beta}\|(1+|X|^2)^{\omega/2}] < \infty$, and for all $x \in \mathcal{X}$, $\frac{\partial \mathcal{E}(x;\beta)}{\partial \beta}$ is continuous around $\beta_0$; (5) $k_n = O(n^{d_x/(2\gamma+d_x)})$, $n^{-\gamma/(2\gamma+d_x)} \times \|\frac{\mathcal{E}(\cdot;\beta_o)}{1-p_o(\cdot)} - \Pi_{2n}\frac{\mathcal{E}(\cdot;\beta_o)}{1-p_o(\cdot)}\|_2 = o(n^{-1/2})$; (6) either (6a) $\sup_{\beta \in B:: \beta-\beta_0: \leq\delta_0} \sup_{x\in\mathcal{X}} \|\mathcal{E}(x,\beta)\| \leq const. < \infty$ for some small $\delta_0 > 0$; or (6b) $E_a[\sup_{\beta \in B: \|\beta-\beta_0\|\leq\delta_0} \|\mathcal{E}(X,\beta)\|^4] \leq const. < \infty$ for some small $\delta_0 > 0$, and $f_{X_a}(\cdot) \in \Lambda_c^\gamma(\mathcal{X})$ with $\gamma > 3d_x/4$; or (6c) $E_a[\sup_{\beta \in B: \|\beta-\beta_0\|\leq\delta_0} \|\mathcal{E}(X,\beta)\|^2] \leq const. < \infty$ for some small $\delta_0 > 0$, and $f_{X_a}(\cdot) \in \Lambda_c^\gamma(\mathcal{X})$ with $\gamma > d_x$.

THEOREM 7. *Let $\widehat{\beta}$ be the IPW-GMM estimator given in* (14) *or* (15). *Under Assumptions* 1, 2, 3 *and* 4, *we have* $\sqrt{n}(\widehat{\beta} - \beta_0) \Rightarrow \mathcal{N}(0,V)$, *with $V$ the same as in Theorem* 4.

REMARK 3. (i) The weighting $\omega$ is needed since the support of the conditioning variable $X$ is assumed to be the entire Euclidean space. When $X$ has bounded support and $f_{X|D=0}$ is bounded above and below over its support, we can simply set $\omega = 0$ in Assumptions 3 and 4 and replace 4.1 with the assumption that 3.1 holds with $\gamma > d_x/2$. Note that Assumption 4.6a is easily satisfied when $X$ has compact support. When $\mathcal{X} = \mathcal{R}^{d_x}$, Assumption 4.6a rules out $\mathcal{E}(x,\beta)$ being linear in $x$; Assumptions 4.6b or 4.6c allow for linear $\mathcal{E}(x,\beta)$ but need smoother propensity score $p(x)$ and density $f_{X|D=0}$. (ii) Assumptions 3 and 4 again allow for non-smooth moment conditions. (iii) Since $\frac{f_{X|D=0}(X)}{f_X(X)} = \frac{1-p_o(X)}{1-p}$, the assumption $0 < \underline{p} \leq p_o(x) \leq \overline{p} < 1$ implies that $\frac{1-\overline{p}}{1-\underline{p}} \leq \frac{f_{X|D=0}(X)}{f_X(X)} \leq \frac{1-\underline{p}}{1-\overline{p}}$, hence $E(\cdot)$ and $E_a(\cdot)$ in Assumptions 3 and 4 are effectively equivalent. (iv) Although Assumption 3.1 imposes the same strong condition $0 < \underline{p} \leq p_o(x) \leq \overline{p} < 1$ as that typically assumed in the program evaluation literature, unlike most existing papers on estimation of average treatment effects, our paper allows for unbounded support of $X$ and assumes weaker smoothness on $p_o(x)$ and $\mathcal{E}(\cdot;\beta_o)$. In particular, if we let $k_n = O(n^{\frac{d_x}{2\gamma+d_x}})$, the growth order which leads to the optimal convergence rate of $\|\widehat{p}(\cdot) - p_o(\cdot)\|_2 = O_p(n^{-\gamma/(2\gamma+d_x)})$, then Assumption 4.5 is satisfied with $\|\frac{\mathcal{E}(\cdot;\beta_o)}{1-p_o(\cdot)} - \Pi_{2n}\frac{\mathcal{E}(\cdot;\beta_o)}{1-p_o(\cdot)}\|_2 = o(n^{-d_x/(2(2\gamma+d_x))}) = o(k_n^{-1/2})$.

4.2. *IPW estimation with parametric or known propensity score.* The case of moment condition (2) is simpler, and therefore, we briefly discuss it first. Theorems 1 and 2 have shown that knowledge about the propensity score does not change the semiparametric efficiency bound. Furthermore, Theorems 4 and 7 show that both a nonparametric CEP-GMM estimator and a nonparametric IPW-GMM estimator for $\beta$ achieve this semiparametric efficiency bound regardless of whether the propensity score is unknown, known or parametrically specified. The following theorem also states, without proof, the interesting result that the parametric IPW estimator using $p(X;\hat{\gamma})$ is in fact less efficient than the one using a nonparametric estimate $\hat{p}(X)$ in (12), but is more efficient than the one using the known $p(X)$.

THEOREM 8. *Suppose that $E[S_\gamma(D,X)S_\gamma(D,X)']$ is positive definite and that the parametric model $p(X_i;\gamma)$ is correctly specified. Under moment condition (2) and using the optimally weighted sample moment condition (12), an IPW-GMM estimator for $\beta$ using a parametric estimate of $p(X_i;\hat{\gamma})$ in place of $\hat{p}(X_i)$ in (12) is more efficient than the one using the known $p(X_i)$, but is less efficient than the one using a nonparametric estimate $\hat{p}(X_i)$ of the propensity score.*

This result is based on the following relations, which hold asymptotically:
$$Avar\left(\frac{\sqrt{n}}{n_a}\sum_{j=1}^{n_a}m(Z_j;\beta_0)\frac{1-\hat{p}}{1-p(X_i;\hat{\gamma})}\right) \leq Avar\left(\frac{\sqrt{n}}{n_a}\sum_{j=1}^{n_a}m(Z_j;\beta_0)\frac{1-\hat{p}}{1-p(X_i)}\right)$$
and
$$Avar\left(\frac{\sqrt{n}}{n_a}\sum_{j=1}^{n_a}m(Z_j;\beta_0)\frac{1-\hat{p}}{1-p(X_i;\hat{\gamma})}\right) \geq Avar\left(\frac{\sqrt{n}}{n_a}\sum_{j=1}^{n_a}m(Z_j;\beta_0)\frac{1-\hat{p}}{1-\hat{p}(X_i)}\right).$$

Now consider the more interesting case where moment condition (1) holds and sample moment condition (11) is used. Consider the case when the parametric propensity score is correctly specified. First, it is clear that the optimally weighted IPW-GMM estimator of $\beta$ based on (11) that uses a nonparametric estimate of $\hat{p}(X)$ does not achieve the efficiency bound in Theorem 3, because we see from Theorem 7 that this estimator achieves instead the variance bound in Theorem 1, which is larger than the variance bound in Theorem 3.

However, the parametric two step IPW estimator that uses a parametric first step for $p(X;\gamma)$ does not achieve the efficiency bound in Theorem 3 either. To see this, note that the parametric two step IPW estimator is based on the moment condition
$$\sqrt{n}\frac{1}{n_a}\sum_{j=1}^{n_a}m(Z_j;\beta)\frac{p(X_j;\hat{\gamma})}{1-p(X_j;\hat{\gamma})}\frac{1-\hat{p}}{\hat{p}},$$

which has a linear influence function representation of
$$\frac{1}{p}\left[m(Z_i;\beta_0)\frac{(1-D_i)p(X_i)}{1-p(X_i)} + \text{Proj}\left(\mathcal{E}(X_i;\beta_0)\frac{D_i-p(X_i)}{1-p(X_i)} \,\Big|\, S_\gamma(D_i,X_i)\right)\right],$$
where
$$\text{Proj}\left(\mathcal{E}(X_i;\beta_0)\frac{D_i-p(X_i)}{1-p(X_i)} \,\Big|\, S_\gamma(D_i,X_i)\right)$$
$$= E\left[\mathcal{E}(X;\beta_0)\frac{p_\gamma(X)}{1-p(X)}\right]$$
$$\times E[S_\gamma(D_i,X_i)S_\gamma(D_i,X_i)']^{-1}S_\gamma(D_i,X_i)$$
is the influence function from the first step estimation of $\gamma$. The difference between this influence function and that in Theorem 3 can be verified to be equal to
$$\text{Res}\left((D-p(X))\frac{p(X)}{1-p(X)}\mathcal{E}(X;\beta_0) \,\Big|\, S_\gamma(D_i,X_i)\right),$$
which is obviously orthogonal to the influence function of Theorem 3. Therefore, the two step parametric IPW estimator has a variance larger than the efficiency bound under the assumption of correct specification of the parametric model for $p(X;\gamma)$.

An IPW type estimator that achieves the efficiency bound under correct specification can be obtained by combining both nonparametric and parametric estimates of the propensity score. Such an efficient moment condition is given by

(16) $$\sqrt{n}\frac{1}{n_a}\sum_{j=1}^{n_a}m(Z_j;\beta)\frac{p(X_j;\hat{\gamma})}{1-\hat{p}(X_j)}\frac{1-\hat{p}}{\hat{p}},$$

where $\hat{\gamma}$ is the maximum likelihood estimator for $\gamma_0$ and $\hat{p}(X)$ is the sieve estimate of the propensity score. This moment condition has the following asymptotic linear representation:
$$\frac{1}{p}\frac{1}{\sqrt{n}}\sum_{i=1}^n\left[(1-D_i)(m(Z_i;\beta_0)-\mathcal{E}(X_i;\beta_0))\frac{p(X_i)}{1-p(X_i)}+p(X_i)\mathcal{E}(X_i;\beta_0)\right]$$
$$+ E\left[\frac{\mathcal{E}(X;\beta_0)}{p}p_\gamma(X_i)\right]\sqrt{n}(\hat{\gamma}-\gamma),$$
which is identical to the influence function under correct parametric specification of $p(X;\gamma)$ leading to the semiparametric efficiency bound in Theorem 3.

The case where the propensity score is fully known can be considered a special case of parametric propensity score where the parameters are known.

In this case, the efficient moment condition is as in (16) after replacing $p(X_j; \hat{\gamma})$ with the known $p(X_j)$.

It is finally worth noting that Assumption 2 is an identification assumption that is not testable. Therefore, both the CEP-GMM estimator and the IPW-GMM estimator will converge to the same population limit regardless of whether Assumption 2 holds, as long as the same weighting matrix is being used. The population difference between CEP and IPW can only arise from the parametric mis-specification of the approximating models for $\mathcal{E}(X; \beta)$ and $p(X)$.

**5. Empirical illustration.** We illustrate our method empirically using data from the Indian National Sample Survey (NSS), which is used to monitor changes in the distribution of private consumption in India. Several researchers have argued that changes in survey methodology caused noncomparability between poverty estimates calculated for 1999–2000 and those from previous years. Changes in the questionnaire likely led to the overestimation of food consumption, and hence to the underestimation of poverty [Deaton and Kozel (2005), Tarozzi (2007)]. In other words, a missing data problem arises because the variable of interest (expenditure recorded using the "standard" questionnaire) is not observed. Deaton and Drèze (2002), Deaton (2003) and Tarozzi (2007) argue that expenditure in a set of miscellaneous items for which the questionnaire was not modified ("comparable items" hereafter) can be used as a proxy variable to produce an estimate of poverty for 1999–2000 that is comparable with previous years.

We assume that the object of interest is the cumulative distribution function (c.d.f.) for rural India in 1999–2000 of a measure of total monthly expenditure that is comparable with previous NSS rounds. In the terminology used in this article, this situation corresponds to a verify-out-of-sample case (1), where the parameter of interest $\beta$ is identified in terms of a variable $Y$ that is not observed in the primary sample (the 1999–2000 survey). The moment function takes the form $m(Z; \beta_0) = 1(Y \leq y) - \beta_0$, where $y$ is a given threshold. We use the previous round of the NSS (1993–94) as auxiliary survey, and expenditure in "comparable items" as proxy variable $X$. The crucial identifying assumption is that the distribution of $Y$ conditional on $X$ remained stable between 1993–94 and 1999–2000 [Tarozzi (2007)].

Table 1 reports point estimates and standard errors for the c.d.f. at selected thresholds. The first column reports estimates using the noncomparable data from the primary sample. Column 2 reports CEP-GMM estimates, calculated using 3rd order polynomial splines in expenditure in comparable items as sieve basis, with 10 knots at the equal range quantiles of the empirical distribution of the proxy variables. Column 3 reports estimates obtained using moment condition (9), but with a nonparametric first step where we estimate $P(X)$ using sieve-logit, including the basis functions we used for

EFFICIENT GMM WITH AUXILIARY DATA    19TABLE 1
*Cumulative distribution functions (×100) of total (log) household expenditure*

| $y$ | (1) Unadjusted (primary) | (2) Adjusted NP CEP | (3) Adjusted NP IPW | (4) Adjusted Par. IPW | (5) Adjusted Eff. CEP |
|---|---|---|---|---|---|
| 6    | 2.92 (0.067)   | 3.388 (0.0695)  | 3.387 (0.0694)  | 3.15 (0.0594)  | 3.23 (0.0598)  |
| 6.25 | 5.67 (0.092)   | 6.521 (0.0948)  | 6.522 (0.0948)  | 6.31 (0.0846)  | 6.38 (0.0845)  |
| 6.50 | 11.06 (0.125)  | 12.272 (0.1237) | 12.273 (0.1234) | 12.21 (0.1165) | 12.21 (0.1149) |
| 6.75 | 20.28 (0.161)  | 21.679 (0.1588) | 21.674 (0.1587) | 21.89 (0.1645) | 21.76 (0.1575) |
| 7    | 34.06 (0.189)  | 35.052 (0.1763) | 35.041 (0.1772) | 35.53 (0.1794) | 35.28 (0.1738) |
| 7.25 | 50.75 (0.200)  | 50.600 (0.1967) | 50.592 (0.1975) | 51.19 (0.1948) | 50.88 (0.1920) |
| 7.50 | 66.98 (0.188)  | 65.682 (0.1925) | 65.687 (0.1929) | 66.15 (0.1973) | 65.91 (0.1880) |

Source: Authors' calculations from Indian National Sample Survey, rounds 50 (1993–94, $n = 58,846$) and 55 (1999–2000, $n = 62,679$), rural sector only from the major Indian states, which account for more than 95% of the total population. Column (1)—Calculated from the unadjusted primary sample. Column (2)—CEP-GMM cubic sieve Estimator, with 10 knots, using "comparable items" as predictor. Column (3)—IPW-GMM. Flexible logit with cubic sieve, with 10 knots, using "comparable items" as predictor. Column (4)—Parametric IPW Estimator. The propensity score is estimated using logit and including total expenditure in "comparable items" as sole predictor. Column (5)—Semiparametric estimator efficient for the case of correctly specified propensity score.

CEP-GMM as regressors. In column 4 we impose a parametric model, and we estimate the propensity score using logit, with $X$ entered linearly in the single index. Column 5 reports the results for the estimator described in Section 3.2, which is efficient when a parametric model is correctly specified for $P(X)$.

For values of $Y$ below 7 the adjusted estimates of the cdf in columns 2 to 4 are slightly larger than the unadjusted figures in column 1. As expected, CEP and IPW non-parametric estimators produce virtually identical results. The estimates in columns 4 and 5 impose a simple logit for the propensity score, but they are still very similar. In the verify-out-of-sample case, knowledge of a parametric form for $P(X)$ lowers the semiparametric efficiency bound, and this may explain why in some cases the standard errors in column 4 are lower than those in columns 2 and 3, where the estimator is only efficient when $P(X)$ is unknown. Note also that when the parametric assumption is correct the efficient estimator is the one in column 5. Indeed the standard errors for this estimator are always lower or virtually identical to those in column 4 every time this latter estimator is more precise than the nonparametric estimators in columns 2 and 3.

**6. Conclusions.** We derive semiparametric efficiency bounds for the estimation of parameters defined through general nonlinear, possibly nonsmooth and over-identified moment conditions, when variables in the primary sample of interest are missing. For identification we rely on the validity of a conditional independence assumption and on the availability of an auxiliary sample that contains information on the relation between missing variables and other proxy variables that are also observed in the primary sample. We study two alternative frameworks. In the first case ("verify-out-of-sample") validation is done with an auxiliary data set which is independent from the primary data set of interest. In the second case ("verify-in-sample") a subset of the observations in the primary sample is validated.

We show that the optimally weighted CEP-GMM estimators achieve the semiparametric efficiency bounds when the propensity score is unknown, or is known or belongs to a correctly specified parametric family. These estimators only use a nonparametric estimate of the conditional expectation of the moment functions, and their asymptotic efficiency is obtained under regularity conditions weaker than the existing ones in the literature. In particular, these CEP-GMM estimators still achieve efficiency bounds when proxy (conditioning) variables have unbounded supports and moment conditions are not smooth.

We also prove that an optimally weighted IPW-GMM estimator is semiparametrically efficient with fully unknown propensity score. However, this estimator is not efficient when the propensity score is either known, or is parametrically estimated using a *correctly* specified parametric model; in such instances, appropriate combinations of nonparametric and parametric estimates of the propensity score are needed to achieve the efficiency bounds.

We have also demonstrated that, from the theoretical point of view, the CEP-GMM estimators are more attractive than the IPW-GMM estimators. Recently and independently Imbens, Newey and Ridder (2005) advocated a similar sieve conditional expectation projection based estimator for the average treatment effect parameter in program evaluation applications. Also, for the estimation of the average treatment effects in missing data models, Wang, Linton and Hardle (2004) suggested that a semiparametrically specified propensity score, such as a single index or a partially linear form, can be used to reduce the curse of dimensionality in the nonparametric estimation of the propensity score. An interesting topic for future research is to study the efficiency implications of these semiparametric restrictions on the propensity score.

## APPENDIX A: CALCULATION OF EFFICIENCY BOUNDS

PROOF OF THEOREM 1. We follow closely the structure of semiparametric efficiency bound derivation of Newey (1990) and Bickel, Klaassen, Ritov and Wellner (1993).

*Case* (1). Consider a parametric path $\theta$ for the joint distribution of $Y, D$ and $X$. Define $p_\theta = P_\theta(D=1)$. The joint density function for $Y, D$ and $X$ is given by

$$(17) \quad f_\theta(y,x,d) = p_\theta^d(1-p_\theta)^{1-d} f_\theta(x \mid D=1)^d f_\theta(x \mid D=0)^{1-d} f(y \mid x)^{1-d}.$$

The resulting score function is given by

$$S_\theta(d,y,x) = \frac{d-p_\theta}{p_\theta(1-p_\theta)} \dot{p}_\theta + (1-d) s_\theta(x \mid D=0)$$
$$+ d s_\theta(x \mid D=1) + (1-d) s_\theta(y \mid x),$$

where $s_\theta(y \mid x) = \frac{\partial}{\partial \theta} \log f_\theta(y \mid x)$, $\dot{p}_\theta = \frac{\partial}{\partial \theta} p_\theta$, $s_\theta(x \mid d) = \frac{\partial}{\partial \theta} \log f_\theta(x \mid d)$. The tangent space of this model is therefore given by:

$$(18) \quad \mathcal{T} = a(d - p_\theta) + (1-d) s_\theta(x \mid D=0) + (1-d) s_\theta(y \mid x) + d s_\theta(x \mid D=1),$$

where $\int s_\theta(y \mid x) f_\theta(y \mid x) \, dy = 0$, $\int s_\theta(x \mid d) f_\theta(x \mid d) \, dx = 0$, and $a$ is a finite constant.

Consider first the case when the model is exactly identified. In this case $\beta$ is uniquely identified by condition (1). Differentiating under the integral gives

$$(19) \quad \frac{\partial \beta(\theta)}{\partial \theta} = -(\mathcal{J}_\beta^1)^{-1} E\left[ m(Z;\beta) \frac{\partial \log f_\theta(Y, X \mid D=1)}{\partial \theta'} \,\bigg|\, D=1 \right].$$

The second component of the right-hand side of this expression can be calculated as

$$(20) \quad E[m(Z;\beta) s_\theta(Y \mid X)' \mid D=1] + E[m(Z;\beta) s_\theta(X \mid D=1)' \mid D=1].$$

Pathwise differentiability follows if we can find $\Psi^1(Y,X,D) \in \mathcal{T}$ such that

$$(21) \quad \partial \beta(\theta)/\partial \theta = E[\Psi^1(Y,X,D) S_\theta(Y,X,D)'].$$

Define $p_\theta = \int p_\theta(x) f_\theta(x) \, dx$, $\mathcal{E}_\theta(X) = E[m(Z;\beta) \mid X]$. It can be verified that pathwise differentiability is satisfied by choosing: $\Psi^1(Y,X,D) = -(\mathcal{J}_\beta^1)^{-1} \times F_\beta^1(Y,X,D)$ where

$$(22) \quad F_\beta^1(Y,X,D) = \frac{1-D}{p} \frac{p(X)}{1-p(X)} [m(Z;\beta) - \mathcal{E}(X)] + \frac{\mathcal{E}(X)}{p} D.$$

Since $\mathcal{J}_\beta^1$ is a nonsingular transformation, this can be shown proving that

$$(23) \quad \begin{aligned} & E\left[ m(Z;\beta) \frac{\partial}{\partial \theta'} \log f_\theta(Y, X \mid D=1) \,\bigg|\, D=1 \right] \\ &= E[F_\beta^1(Y,X,D) S_\theta(Y,X,D)']. \end{aligned}$$

This can in turn be verified by checking that

$$E[m(Z;\beta)s_\theta(Y \mid X)' \mid D = 1]$$
$$= E\left[\frac{1-D}{p}\frac{p(X)}{1-p(X)}[m(Z;\beta) - \mathcal{E}(X)]s_\theta(Y \mid X)'\right],$$
$$E[m(Z;\beta)s_\theta(X \mid D = 1)' \mid D = 1] = E\left[\frac{\mathcal{E}(X)}{p}Ds_\theta(X \mid D = 1)'\right].$$

Now one can also verify that $F_\beta^1(Y, X, D)$ belongs to the tangent space $\mathcal{T}$ in equation (18), with the first and second terms of $F_\beta^1(Y, X, D)$ taking the role of $(1-d)s_\theta(y|x)$ and $ds_\theta(X|D=1)$, respectively, and the two other components in (18) being identically equal to 0.

Therefore all the conditions of Theorem 3.1 in Newey (1990) hold, so that $\Psi^1$ is the efficient score function and the efficiency bound for regular estimators of the parameter $\beta$ is given by

(24) $\quad V_1 = (\mathcal{J}_\beta^1)^{-1} E[F_\beta^1(Y, X, D)F_\beta^1(Y, X, D)'](\mathcal{J}_\beta^1)'^{-1} = (\mathcal{J}_\beta^1)^{-1}\Omega_\beta^1(\mathcal{J}_\beta^1)'^{-1}.$

*Case* (2). For this case we use an alternative factorization of the likelihood function. Define $p_\theta(x) = P_\theta(D = 1 \mid x)$. The joint density function for $Y$, $D$ and $X$ is given by

(25) $\quad f_\theta(y, x, d) = f_\theta(x)p_\theta(x)^d[1 - p_\theta(x)]^{1-d}f_\theta(y \mid x)^{1-d}.$

The resulting score function is then given by

$$S_\theta(d, y, x) = (1-d)s_\theta(y \mid x) + \frac{d - p_\theta(x)}{p_\theta(x)(1 - p_\theta(x))}\dot{p}_\theta(x) + t_\theta(x),$$

where

$$s_\theta(y \mid x) = \frac{\partial}{\partial\theta}\log f_\theta(y \mid x), \qquad \dot{p}_\theta(x) = \frac{\partial}{\partial\theta}p_\theta(x), \qquad t_\theta(x) = \frac{\partial}{\partial\theta}\log f_\theta(x).$$

The tangent space of this model is therefore given by:

(26) $\quad \mathcal{T} = \{(1-d)s_\theta(y \mid x) + a(x)(d - p_\theta(x)) + t_\theta(x)\}$

where $\int s_\theta(y \mid x)f_\theta(y \mid x)\,dy = 0$, $\int t_\theta(x)f_\theta(x)\,dx = 0$, and $a(x)$ is any square integrable function.

In case (2), equation (19) is replaced by

(27) $\quad \begin{aligned}\frac{\partial\beta(\theta)}{\partial\theta} &= -(\mathcal{J}_\beta^2)^{-1}E\left[m(Z;\beta)\frac{\partial\log f_\theta(Y, X)}{\partial\theta'}\right] \\ &= -(\mathcal{J}_\beta^2)^{-1}\{E[m(Z;\beta)s_\theta(Y \mid X)'] + E[\mathcal{E}(X)t_\theta(X)']\}.\end{aligned}$

Now we replace $F_\beta^1(Y, X, D)$ in (22) with the following:

(28) $\quad F_\beta^2(Y, X, D) = \frac{1-D}{1-p(X)}[m(Z;\beta) - \mathcal{E}(X)] + \mathcal{E}(X)$

and then it can be verified that $E[F_\beta^2(Y,X,D)S_\theta(Y,X,D)'] = E[m(Z;\beta) \times \frac{\partial \log f_\theta(Y,X)}{\partial \theta'}]$. Then the efficient influence function for case (2) is equal to $-(\mathcal{J}_\beta^2)^{-1} F_\beta^2(Y,X,D)$ with the two terms being orthogonal to each other, and the second result in Theorem 1 follows.

Now consider overidentified moment conditions. We only consider case (1), as the derivation for case (2) is analogous. When $d_m > d_\beta$, the moment conditions in (1) is equivalent to the requirement that for any matrix $\mathcal{A}$ of dimension $d_\beta \times d_m$ the following exactly identified system of moment conditions holds $\mathcal{A}E[m(Z;\beta) \mid D=1] = 0$. Differentiating again,

$$\frac{\partial \beta(\theta)}{\partial \theta} = -\left(\mathcal{A}E\left[\frac{\partial m(Z;\beta)}{\partial \beta}\Big|D=1\right]\right)^{-1}$$
$$\times E\left[\mathcal{A}m(Z;\beta)\frac{\partial \log f_\theta(Y,X \mid D=1)}{\partial \theta'}\,\Big|\, D=1\right].$$

Therefore, any regular estimator for $\beta$ is asymptotically linear with influence function of the form

$$-\left(\mathcal{A}E\left[\frac{\partial m(Z;\beta)}{\partial \beta}\Big|D=1\right]\right)^{-1}\mathcal{A}m(z;\beta).$$

For a given matrix $\mathcal{A}$, the projection of the above influence function onto the tangent set follows from the previous calculations, and is given by $-[\mathcal{A}\mathcal{J}_\beta^1]^{-1}F_\beta^1(y,x,d)$. The asymptotic variance corresponding to this efficient influence function for fixed $\mathcal{A}$ is therefore

$$[\mathcal{A}\mathcal{J}_\beta^1]^{-1}\mathcal{A}\Omega\mathcal{A}'[\mathcal{J}_\beta^{1\prime}\mathcal{A}']^{-1}, \tag{29}$$

where $\Omega = E[F_\beta^1(Y,X,D)F_\beta^1(Y,X,D)']$ as calculated above. Therefore, the efficient influence function is obtained when $\mathcal{A}$ minimizes (29). It is easy to show that such matrix $\mathcal{A}$ is equal to $\mathcal{J}_\beta^{1\prime}\Omega^{-1}$, so that the asymptotic variance becomes $V = (\mathcal{J}_\beta^{1\prime}\Omega^{-1}\mathcal{J}_\beta^1)^{-1}$. In fact, a standard textbook calculation shows

$$\mathcal{J}_\beta^{1\prime}\Omega^{-1}\mathcal{J}_\beta^1 - \mathcal{J}_\beta'\mathcal{A}'(\mathcal{A}\Omega\mathcal{A}')^{-1}\mathcal{A}\mathcal{J}_\beta$$
$$= (\mathcal{J}_\beta^{1\prime}\Omega^{-1/2} - \mathcal{J}_\beta^{1\prime}\Omega^{-1/2}\Omega^{1/2\prime}(\Omega^{1/2}\Omega^{1/2\prime})^{-1}\Omega^{1/2})$$
$$\times (\Omega^{-1/2}\mathcal{J}_\beta^1 - \Omega^{1/2\prime}[\Omega^{1/2}\Omega^{1/2\prime}]^{-1}\Omega^{1/2}\Omega^{-1/2}\mathcal{J}_\beta^1) \geq 0. \qquad \square$$

PROOF OF THEOREM 2. As for Theorem 1, it suffices to present the proof for the case of exact identification, since the overidentified case follows from choosing the optimal linear combination matrix. If the propensity score $p(x)$ is known, the score becomes [cf. Hahn (1998)] $S_\theta(d,y,x) = (1-d)s_\theta(y \mid x) + t_\theta(x)$, so that the tangent space becomes $\mathcal{T} = \{(1-d)s_\theta(y \mid x) + t_\theta(x)\}$ where $\int s_\theta(y \mid x) f_\theta(y \mid x) \, dy = 0$, and $\int t_\theta(x) f_\theta(x) \, dx = 0$. Consider case (1) first. The pathwise derivative becomes

$$E\left[\frac{p(X)}{p} m(Z;\beta) s(Y \mid X)'\right] + E\left[\frac{p(X)}{p} \mathcal{E}(X) t(X)'\right].$$

Pathwise differentiability is established by verifying that equation (21) holds, with

(30) $$F_\beta^1(y, x, d) = \frac{1-d}{p} \frac{p(x)}{1-p(x)} (m(z;\beta) - \mathcal{E}(x)) + \frac{\mathcal{E}(x)}{p} p(x).$$

Then the efficient influence function is as before equal to $-(\mathcal{J}_\beta^1)^{-1} F_\beta^1(y, x, d)$, and the result of Theorem 2 follows using Theorem 3.1 of Newey (1990).

Since $p(x)$ does not enter the definition of $\beta$ in case (2), there is no change to the efficient influence function and to the semiparametric efficiency bound for that case. □

PROOF OF THEOREM 3. When $p(X)$ belongs to a correctly specified parametric family $p(X;\gamma)$, the score function for moment (1) becomes

$$S_\theta(d, y, x) = (1-d) s_\theta(y \mid x) + \frac{d - p_\theta(x)}{p_\theta(x)(1 - p_\theta(x))} \frac{\partial p(x;\gamma)}{\partial \gamma'} \frac{\partial \gamma}{\partial \theta} + t_\theta(x).$$

The tangent space is therefore $\mathcal{T} = \{(1-d) s_\theta(y \mid x) + c' S_\gamma(d; x) + t_\theta(x)\}$ where $c$ is a finite vector of constants and $S_\gamma(d; x)$ is the parametric score function. Now define $F_\beta^1(Y, X, D)$ as

$$\frac{1-D}{p} \frac{p(X)}{1-p(X)} [m(Z;\beta) - \mathcal{E}(X)] + \mathrm{Proj}\Big(\mathcal{E}(X) \frac{D - p(X)}{p} \Big| S_\gamma(D, X)\Big).$$

It is clear that $F_\beta^1(Y, X, D)$ lies in the tangent space. Also note that $\frac{\partial \beta(\theta)}{\partial \theta}$ can be written as

$$-(\mathcal{J}_\beta^1)^{-1}\Big\{E[m(Z;\beta) s_\theta(Y \mid X)' \mid D = 1]$$
$$+ E\Big[m(Z;\beta)\Big(t_\theta(x)' + S_\gamma(d;x)' \frac{\partial \gamma}{\partial \theta}\Big) \Big| D = 1\Big]\Big\}.$$

The second term in curly brackets can also be written as

$$\frac{E(D - p(X)) \mathcal{E}(X) S_\gamma(D; X)'}{p} \frac{\partial \gamma}{\partial \theta} + \frac{p(X) \mathcal{E}(X) t_\theta(X)}{p}.$$

With these calculations it can be verified that

$$\frac{\partial \beta(\theta)}{\partial \theta} = -(\mathcal{J}_\beta^1)^{-1} E[F_\beta^1(Y, X, D) S_\theta(Y, X, D)].$$

In particular,

$$E\left[\frac{(D-p(X))\mathcal{E}(X)S_\gamma(D;X)'}{p}\right]$$
$$= E\left[\text{Proj}\left(\mathcal{E}(X)\frac{D-p(X)}{p} \mid S_\gamma(D,X)\right)S_\theta(Y,X,D)'\right].$$

Therefore $-(\mathcal{J}_\beta^1)^{-1}F_\beta^1(Y,X,D)$ is the desired efficient influence function and its variance is given as the efficient variance of Theorem 3. $\square$

## APPENDIX B: PROOFS OF ASYMPTOTIC PROPERTIES

In this Appendix we establish the large sample properties for the IPW-GMM estimator with nonparametrically estimated propensity score function. Again to stress the fact that the true propensity score is unknown, in this Appendix we denote the true propensity score by $p_o(x) \equiv E[D|X=x]$ and any candidate function by $p(x)$.

Denote

$$\mathcal{L}_2(\mathcal{X}) = \left\{h: \mathcal{X} \to \mathcal{R}: \|h\|_2 = \sqrt{\int h(x)^2 f_X(x)\,dx} < \infty\right\}$$

and

$$\mathcal{L}_{2,a}(\mathcal{X}) = \left\{h: \mathcal{X} \to \mathcal{R}: \|h\|_{2,a} = \sqrt{\int h(x)^2 f_{X_a}(x)\,dx} < \infty\right\}$$

as the two Hilbert spaces. We use $\|h\|_2 \asymp \|h\|_{2,a}$ to mean that there are two positive constants $c_1, c_2$ such that $c_1\|h\|_2 \leq \|h\|_{2,a} \leq c_2\|h\|_2$, which is true under the assumption $0 < \underline{p} \leq p_o(x) \leq \overline{p} < 1$.

Proposition B.1 provides large sample properties for the sieve LS estimator $\widehat{p}(x)$ of $p_o(x)$.

PROPOSITION B.1. *Under Assumptions* 3.1, 3.2 *and* 3.5, *and* $\frac{k_n}{n} \to 0$, $k_n \to \infty$, *we have* (i)

$$\|\widehat{p}(\cdot) - p_o(\cdot)\|_{\infty,\omega} = o_p(1); \qquad \|\widehat{p}(\cdot) - p_o(\cdot)\|_{2,a} \asymp \|\widehat{p}(\cdot) - p_o(\cdot)\|_2 = o_p(1);$$

(ii) *in addition, if Assumption* 4.1 *holds, then*

$$\|\widehat{p}(\cdot) - p_o(\cdot)\|_{2,a} \asymp \|\widehat{p}(\cdot) - p_o(\cdot)\|_2 = O_p\left(\sqrt{\frac{k_n}{n}} + (k_n)^{-\gamma/d_x}\right).$$

PROOF. (i) Recall that $\widehat{p}(x)$ is the sieve LS estimator of $p_o(\cdot) \in \Lambda_c^\gamma(\mathcal{X})$ based on the entire sample. That is,

$$\widehat{p}(\cdot) = \arg\min_{p(\cdot) \in \mathcal{H}_n} \frac{1}{n}\sum_{i=1}^n \{D_i - p(X_i)\}^2/2,$$

where $\mathcal{H}_n$ increases with sample size $n$, and is dense in $\Lambda_c^\gamma(\mathcal{X})$ as $k_n \to \infty$ (by Assumption 3.5). Moreover, by Assumptions 3.1 and 3.2 we have the following results: (1) the parameter space is compact under the norm $\|\cdot\|_{\infty,\omega}$ for $\omega > 0$, see Ai and Chen (2003); (2) $E[\{D_i - p(X_i)\}^2/2]$ is uniquely maximized at $p_o(x) = E[D|X = x] \in \mathcal{H}$; (3) $E[\{D_i - p(X_i)\}^2/2]$ is continuous in $p(\cdot)$ under the metric $\|\cdot\|_{\infty,\omega}$; and (4)

$$\sup_{p(\cdot)\in\mathcal{H}}\left|\frac{1}{n}\sum_{i=1}^n \{D_i - p(X_i)\}^2/2 - E\{D_i - p(X_i)\}^2/2\right| = o_p(1);$$

where both results (3) and (4) are due to the fact that for any $p(\cdot), \widetilde{p}(\cdot) \in \mathcal{H}$,

$$|\{D_i - p(X_i)\}^2 - \{D_i - \widetilde{p}(X_i)\}^2|$$
$$= |\{2D_i - [p(X_i) + \widetilde{p}(X_i)]\}[p(X_i) - \widetilde{p}(X_i)]|$$
$$\le const. |[p(X_i) - \widetilde{p}(X_i)](1 + X_i'X_i)^{-\omega/2}| \times (1 + X_i'X_i)^{\omega/2}.$$

Now $E[(1 + X_i'X_i)^{\omega/2}] < \infty$ by Assumption 3.2.

Hence by either Theorem 0 in Gallant and Nychka (1987) or Lemma 2.9 and Theorem 2.1 in Newey (1994), $\|\widehat{p}(\cdot) - p_o(\cdot)\|_{\infty,\omega} = o_p(1)$. Now

$$\|\widehat{p}(\cdot) - p_o(\cdot)\|_2 = \sqrt{\int [\widehat{p}(x) - p_o(x)]^2 f_X(x)\,dx}$$
$$\le \sqrt{(\|\widehat{p}(\cdot) - p_o(\cdot)\|_{\infty,\omega})^2 \int (1 + x'x)^\omega f_X(x)\,dx} = o_p(1)$$

(by Assumption 3.2).

(ii) We can obtain the convergence rate of $\|\widehat{p}(\cdot) - p_o(\cdot)\|_2$ by applying Theorem 1 in Chen and Shen (1998) or Theorem 2 in Shen and Wong (1994). Let $L_n(p(\cdot)) = \frac{1}{n}\sum_{i=1}^n \ell(D_i, X_i, p(\cdot))$ with $\ell(D_i, X_i, p(\cdot)) = -\{D_i - p(X_i)\}^2/2$. Since all the assumptions of Chen and Shen (1998) Theorem 1 are satisfied given our Assumptions 3.1 and 3.2. We obtain

$$\|\widehat{p}(\cdot) - p_o(\cdot)\|_2 = O_p\left(\max\left\{\sqrt{\frac{k_n}{n}}, \|p_o - \Pi_{2n}p_o\|_2\right\}\right).$$

Under Assumption 4.1, for $p_o \in \Lambda_c^\gamma(\mathcal{X})$, there exists $\Pi_{\infty n}p_o \in \Lambda_c^\gamma(\mathcal{X})$ such that for any fixed $\omega > \gamma$,

$$\|p_o - \Pi_{\infty n}p_o\|_{\infty,\omega} = \sup_x |[p_o(x) - \Pi_{\infty n}p_o(x)](1 + |x|^2)^{-\omega/2}|$$
$$\le const. (k_n)^{-\gamma/d_x},$$

see Ai and Chen (2003). Hence by Assumption 4.1 with $\omega = \gamma + \epsilon$ for a small $\epsilon > 0$,

$$\|p_o - \Pi_{2n}p_o\|_2 \leq \|p_o - \Pi_{\infty n}p_o\|_2$$
$$= \sqrt{\int [p_o(x) - \Pi_{\infty n}p_o(x)]^2 f_X(x)\,dx}$$
$$\leq \sqrt{(\|p_o(\cdot) - \Pi_{\infty n}p_o(\cdot)\|_{\infty,\omega})^2 \int (1+x'x)^\omega f_X(x)\,dx}$$
$$\leq c'(k_n)^{-\gamma/d_x}.$$

Then $\|\widehat{p}(\cdot) - p_o(\cdot)\|_2 = O_p(\sqrt{\frac{k_n}{n}} + (k_n)^{-\gamma/d_x}) = o_p(1)$. □

PROOF OF THEOREM 6. We only provide the proof of the IPW-GMM estimator for moment condition (1), since the one for moment condition (2) is very similar. We establish this theorem by applying Theorem 1 in Chen, Linton and van Keilegom (2003) (hereafter CLK) with their $\theta$ being our $\beta$ and their $h$ being our $p(\cdot)$. Define

$$M_n(\beta, p(\cdot)) = \frac{1}{n_a} \sum_{i=1}^{n_a} m(Z_i, \beta) \frac{p(X_i)}{1 - p(X_i)};$$

$$M(\beta, p(\cdot)) = E_a\left[m(Z_i, \beta) \frac{p(X_i)}{1 - p(X_i)}\right] = E\left[m(Z, \beta) \frac{p(X)}{1 - p(X)} \Big| D = 0\right].$$

CLK's conditions (1.1) and (1.2) are directly implied by our Assumptions 1.1, 2 and moment condition (1). Note that for any $p(\cdot) \in \mathcal{H}$, $0 < \frac{1}{1-\underline{p}} \leq \frac{1}{1-p(X)} \leq \frac{1}{1-\overline{p}} < \infty$, we have

$$|M(\beta, p(\cdot)) - M(\beta, p_o(\cdot))|$$
$$= \left|E\left[m(Z, \beta)\left\{\frac{p(X)}{1-p(X)} - \frac{p_o(X)}{1-p_o(X)}\right\} \Big| D = 0\right]\right|$$
$$\leq \frac{1}{(1-\overline{p})^2} E_a[\|m(Z,\beta)\|(1+|X|^2)^{\omega/2}]$$
$$\quad \times \sup_{x \in \mathcal{X}} |[p(x) - p_o(x)](1+|x|^2)^{-\omega/2}|$$
$$\leq \frac{1}{(1-\overline{p})^2} \left\{E_a\left[\sup_{\beta \in B}\|m(Z,\beta)\|^2\right] \times E_a[(1+|X|^2)^\omega]\right\}^{1/2}$$
$$\quad \times \|p(\cdot) - p_o(\cdot)\|_{\infty,\omega},$$

where the last inequality is due to our Assumptions 3.1, 3.2 and 3.4, hence CLK's condition (1.3) is satisfied with respect to the norm $\|\cdot\|_{\mathcal{H}} = \|\cdot\|_{\infty,\omega}$.

CLK's condition (1.4) $\|\widehat{p}(\cdot) - p_o(\cdot)\|_{\infty,\omega} = o_p(1)$ is implied by Proposition B.1(i). Note that

$$E_a\left[\sup_{\|\beta-\widetilde{\beta}\|<\delta, \|p(\cdot)-\widetilde{p}(\cdot)\|_{\infty,\omega}<\delta}\left|m(Z_i,\beta)\frac{p(X_i)}{1-p(X_i)} - m(Z_i,\widetilde{\beta})\frac{\widetilde{p}(X_i)}{1-\widetilde{p}(X_i)}\right|\right]$$

$$\leq E_a\left[\sup_{\|\beta-\widetilde{\beta}\|<\delta}\|m(Z_i,\beta) - m(Z_i,\widetilde{\beta})\| \times \sup_{p(\cdot)\in\mathcal{H}}\left|\frac{p(X_i)}{1-p(X_i)}\right|\right]$$

$$+ E_a\left[\sup_{\widetilde{\beta}\in B}\|m(Z_i,\widetilde{\beta})\| \times \sup_{\|p(\cdot)-\widetilde{p}(\cdot)\|_{\infty,\omega}<\delta}\left|\frac{p(X_i)}{1-p(X_i)} - \frac{\widetilde{p}(X_i)}{1-\widetilde{p}(X_i)}\right|\right]$$

$$\leq E_a\left[\sup_{\|\beta-\widetilde{\beta}\|<\delta}\|m(Z_i,\beta) - m(Z_i,\widetilde{\beta})\|\right] \times \frac{\overline{p}}{1-\overline{p}}$$

$$+ E_a\left[\sup_{\widetilde{\beta}\in B}\|m(Z_i,\widetilde{\beta})\|(1+|X_i|^2)^{\omega/2}\right]$$

$$\times \frac{\sup_{\|p(\cdot)-\widetilde{p}(\cdot)\|_{\infty,\omega}<\delta}\sup_{x\in\mathcal{X}}|[p(x)-\widetilde{p}(x)](1+|x|^2)^{-\omega/2}|}{(1-\overline{p})^2}$$

$$\leq const.b(\delta) + const.\delta,$$

where the last inequality is due to our Assumptions 3.1–3.4 and Proposition B.1(i). Then CLK's condition (1.5) is satisfied, hence $\widehat{\beta} - \beta_0 = o_p(1)$. □

LEMMA B.2. *Under Assumptions* 1, 2, 3 *and* 4, *we have*

$$\sqrt{n}E\left\{\mathcal{E}(X,\beta_o)\frac{\widehat{p}(X) - p_o(X)}{1 - p_o(X)}\right\}$$

$$= \frac{1}{\sqrt{n}}\sum_{i=1}^{n}\frac{D_i - p_o(X_i)}{1 - p_o(X_i)}\mathcal{E}(X,\beta_o) + o_p(1).$$

PROOF. To establish this result, we follow the approach in Shen (1997) and Chen and Shen (1998). Recall $p_o(x) = E[D|X=x] \in \Lambda_c^\gamma(\mathcal{X})$ and

$$\widehat{p}(\cdot) = \arg\min_{p(\cdot)\in\mathcal{H}_n} \frac{1}{n}\sum_{i=1}^{n}\{D_i - p(X_i)\}^2/2.$$

Define the inner product associated with the space $\mathcal{L}_2(\mathcal{X})$ as

$$\langle h, g \rangle = E\{h(X)g(X)\} \qquad \text{hence } \|h(\cdot)\|_2^2 = \langle h, h \rangle = E[\{h(X)\}^2].$$

Then the Riesz representor $v^*$ for functional $E\{\mathcal{E}(X,\beta_o)\frac{p(X)-p_o(X)}{1-p_o(X)}\}$ is simply given by

$$v^*(X) = \frac{\mathcal{E}(X,\beta_o)}{1 - p_o(X)},$$

this is because
$$\|v^*\|^2 = \sup_{p(\cdot) \in \mathcal{H}: p \neq p_o} \frac{[E\{\mathcal{E}(X, \beta_o)((p(X) - p_o(X))/(1 - p_o(X)))\}]^2}{E[(p(X) - p_o(X))^2]}$$
$$= E\left[\left(\frac{\mathcal{E}(X, \beta_o)}{1 - p_o(X)}\right)^2\right]$$

and
$$E\left\{\mathcal{E}(X, \beta_o)\frac{p(X) - p_o(X)}{1 - p_o(X)}\right\} = \langle v^*, p(\cdot) - p_o(\cdot) \rangle$$
$$= E\{v^*(X)[p(X) - p_o(X)]\}.$$

Let $L_n(p(\cdot)) = \frac{1}{n}\sum_{i=1}^n \ell(D_i, X_i, p(\cdot))$ with $\ell(D_i, X_i, p(\cdot)) = -\{D_i - p(X_i)\}^2/2$. Let $U_i \equiv D_i - p_o(X_i)$. Then by definition $E[U_i|X_i] = 0$, and $\ell(D_i, X_i, p(\cdot)) = -\{U_i - [p(X_i) - p_o(X_i)]\}^2/2$. We denote $\mu_n(g) = \frac{1}{n} \times \sum_{i=1}^n [g(D_i, X_i) - E(g(D_i, X_i))]$ as the empirical process indexed by $g$, and $\varepsilon_n$ be any positive sequence with $\varepsilon_n = o(\frac{1}{\sqrt{n}})$. Then by definition,

$$0 \leq L_n(\widehat{p}) - L_n(\widehat{p} \pm \varepsilon_n \Pi_{2n} v^*)$$
$$= \mu_n(\ell(D_i, X_i, \widehat{p}) - \ell(D_i, X_i, \widehat{p} \pm \varepsilon_n \Pi_{2n} v^*))$$
$$+ E(\ell(D_i, X_i, \widehat{p}) - \ell(D_i, X_i, \widehat{p} \pm \varepsilon_n \Pi_{2n} v^*)).$$

A simple calculation yields
$$E(\ell(D_i, X_i, \widehat{p}) - \ell(D_i, X_i, \widehat{p} \pm \varepsilon_n \Pi_{2n} v^*))$$
$$= \pm \varepsilon_n E[\Pi_{2n} v^*(X_i)\{\widehat{p}(X_i) - p_o(X_i)\}]$$
$$+ \tfrac{1}{2}\varepsilon_n^2 E[\{\Pi_{2n} v^*(X_i)\}^2],$$
$$\mu_n(\ell(D_i, X_i, \widehat{p}) - \ell(D_i, X_i, \widehat{p} \pm \varepsilon_n \Pi_{2n} v^*))$$
$$= \mp \varepsilon_n \times \mu_n(\Pi_{2n} v^* U_i)$$
$$\pm \varepsilon_n \times \mu_n\left(\Pi_{2n} v^* \frac{2\{\widehat{p}(\cdot) - p_o(\cdot)\} \pm \varepsilon_n \Pi_{2n} v^*}{2}\right)$$

hence
$$0 \leq \mp \mu_n(\Pi_{2n} v^*(X_i) U_i) \pm E[\Pi_{2n} v^*(X_i)\{\widehat{p}(X_i) - p_o(X_i)\}]$$
$$\pm \mu_n(\Pi_{2n} v^*(X_i)\{\widehat{p}(X_i) - p_o(X_i)\}) + \frac{\varepsilon_n}{2n} \sum_{i=1}^n \{\Pi_{2n} v^*(X_i)\}^2$$
$$= \mp \mu_n([\Pi_{2n} v^* - v^*] U_i) \pm \mu_n(v^* U_i)$$
$$\pm E[[\Pi_{2n} v^* - v^*]\{\widehat{p} - p_o\}] \mp E[v^*\{\widehat{p} - p_o\}]$$
$$\pm \mu_n(\Pi_{2n} v^*(X_i)\{\widehat{p}(X_i) - p_o(X_i)\}) + \frac{\varepsilon_n}{2n} \sum_{i=1}^n \{\Pi_{2n} v^*(X_i)\}^2.$$

In the following we shall establish (B2.1)–(B2.4):

(B2.1) $$\mu_n([\Pi_{2n}v^*(X_i) - v^*(X_i)]U_i) = o_p\left(\frac{1}{\sqrt{n}}\right),$$

(B2.2) $$E([\Pi_{2n}v^*(X_i) - v^*(X_i)]\{\widehat{p}(X_i) - p_o(X_i)\}) = o_p\left(\frac{1}{\sqrt{n}}\right),$$

(B2.3) $$\mu_n(\Pi_{2n}v^*(X_i)\{\widehat{p}(X_i) - p_o(X_i)\}) = o_p\left(\frac{1}{\sqrt{n}}\right),$$

(B2.4) $$\frac{1}{n}\sum_{i=1}^n \{\Pi_{2n}v^*(X_i)\}^2 = O_p(1).$$

Note that (B2.1) is implied by Chebychev inequality, i.i.d. data, and $\|\Pi_n v^* - v^*\|_2 = o(1)$ which is satisfied given the expression for $v^*$ and Assumptions 3.1 and 4.5. (B2.2) is implied by Assumption 4.5 and $\|\widehat{p}(\cdot) - p_o(\cdot)\|_2 = O_p(n^{-\gamma/(2\gamma+d_x)})$ from Proposition B.1(ii). (B2.4) is implied by Markov inequality, i.i.d. data, and Assumptions 3.1 and 4.5. Finally for (B2.3), let $\mathcal{F}_n = \{\Pi_{2n}v^*(\cdot)h(\cdot) : h(\cdot) \in \Lambda_c^\gamma(\mathcal{X})\}$, then by Assumption 4.1, $\log N_{[\cdot]}(\delta, \mathcal{F}_n, \|\cdot\|_2) \leq const.(\frac{c}{\delta})^{d_x/\gamma}$ for any $\delta > 0$. Applying Theorem 3 in Chen and Shen (1998) with their $\delta_n = n^{-\gamma/(2\gamma+d_x)}$, we have

$$\sup_{h \in \mathcal{F}_n : \|h(\cdot) - p_o(\cdot)\|_2 \leq \delta_n} |\sqrt{n}\mu_n(\Pi_{2n}v^*\{h(\cdot) - p_o(\cdot)\})|$$
$$= O_p(n^{-(2\gamma-d_x)/(2(2\gamma+d_x))}) = o_p(1).$$

Hence we obtain (B2.3). Now (B2.1)–(B2.4) imply $0 \leq \pm\mu_n(v^*U_i) \mp E[v^*\{\widehat{p} - p_o\}] + o_p(\frac{1}{\sqrt{n}})$, that is $\sqrt{n}E[v^*(X)\{\widehat{p}(X) - p_o(X)\}] = \frac{1}{\sqrt{n}}\sum_{i=1}^n v^*(X_i)U_i + o_p(1)$, hence the result follows. $\square$

PROOF OF THEOREM 7. Again we only provide the proof of the IPW-GMM estimator for moment condition (1). We establish this theorem by applying Theorem 2 in CLK (2003). Given the definition of $\beta_0$ and Theorem 6, CLK's condition (2.1) is directly satisfied. Note that their $\Gamma_1(\beta, p_o) = \frac{p}{1-p}J_\beta^1$, hence their condition (2.2) is satisfied with our Assumption 1.1.

Following the proof of CLK's Theorem 2, we note that the conclusion of CLK's Theorem 2 remains true when CLK's conditions (2.3)(i) and (2.4) are replaced by the following one:

(*) $$\sup_{\beta \in B\|\beta-\beta_0\|\leq\delta_0} \|M(\beta, \widehat{p}(\cdot)) - M(\beta, p_o(\cdot)) - \Gamma_2(\beta, p_o)[\widehat{p}(\cdot) - p_o(\cdot)]\| = o_p(n^{-1/2}),$$

where

$$\Gamma_2(\beta, p_o)[p(\cdot) - p_o(\cdot)] = E\left\{m(Z,\beta)\frac{p(X) - p_o(X)}{(1 - p_o(X))^2}\Big| D = 0\right\}$$

$$= E_a\left\{\mathcal{E}(X,\beta)\frac{p(X) - p_o(X)}{(1 - p_o(X))^2}\right\}$$

$$= E\left\{\mathcal{E}(X,\beta)\frac{p(X) - p_o(X)}{(1 - p_o(X))^2}\frac{f_{X|D=0}(X)}{f_X(X)}\right\}$$

$$= \frac{1}{1-p}E\left\{\mathcal{E}(X,\beta)\frac{p(X) - p_o(X)}{1 - p_o(X)}\right\},$$

and the last equality is due to $f_{X|D=0}(X)/f_X(X) = (1 - p_o(X))/(1 - p)$.

Before we apply Assumptions 4.6a or 4.6b or 4.6c to verify condition (*), let us check CLK's conditions (2.3)(ii), (2.5) and (2.6). Since for all $\beta$ with $\|\beta - \beta_0\| \leq \delta_0$ and all $p(\cdot)$ with $\|p(\cdot) - p_o(\cdot)\|_{\infty,\omega} \leq \delta_0$, we have

$$|\Gamma_2(\beta, p_o)[p(\cdot) - p_o(\cdot)] - \Gamma_2(\beta_o, p_o)[p(\cdot) - p_o(\cdot)]|$$

$$= \left|\frac{1}{1-p}E\left\{[\mathcal{E}(X,\beta) - \mathcal{E}(X,\beta_0)]\frac{p(X) - p_o(X)}{1 - p_o(X)}\right\}\right|$$

$$= \left|\frac{\beta - \beta_0}{1-p}E\left\{\frac{\partial \mathcal{E}(X,\overline{\beta})}{\partial \beta}\frac{p(X) - p_o(X)}{1 - p_o(X)}\right\}\right|$$

$$\leq \frac{\|\beta - \beta_0\|}{(1-p)(1-\overline{p})}E\left[\left\|\frac{\partial \mathcal{E}(X,\overline{\beta})}{\partial \beta}\right\|(1 + |X|^2)^{\omega/2}\right]$$

$$\times \sup_{x \in \mathcal{X}}|[p(x) - p_o(x)](1 + |x|^2)^{-\omega/2}|,$$

where $\overline{\beta}$ is in between $\beta$ and $\beta_0$. Thus, under our Assumptions 3.2, 4.4, Proposition B.1(i) and Theorem 6, $|\Gamma_2(\beta, p_o)[p(\cdot) - p_o(\cdot)] - \Gamma_2(\beta_o, p_o)[p(\cdot) - p_o(\cdot)]| \leq const.\|\beta - \beta_0\| \times \|p(\cdot) - p_o(\cdot)\|_{\infty,\omega}$ hence CLK's condition (2.3)(ii) is satisfied.

Now we verify CKL's condition (2.5) by applying their Theorem 3. In fact, given our Theorem 6 and Proposition B.1(i), it suffices to consider some neighborhood around $(\beta_o, p_o)$. Let $\delta_0 > 0$ be a small value, then for all $(\widetilde{\beta}, \widetilde{p}) \in B \times \mathcal{H}$ with $\|\widetilde{\beta} - \beta_o\| \leq \delta_0$ and $\|\widetilde{p} - p_o\|_{\infty,\omega} \leq \delta_0$, we have for any $\delta \in (0, \delta_0]$,

$$E_a\left[\sup_{\|\beta - \widetilde{\beta}\| < \delta, \|p(\cdot) - \widetilde{p}(\cdot)\|_{\infty,\omega} < \delta}\left|m(Z_i, \beta)\frac{p(X_i)}{1 - p(X_i)} - m(Z_i, \widetilde{\beta})\frac{\widetilde{p}(X_i)}{1 - \widetilde{p}(X_i)}\right|^2\right]$$

$$\leq E_a\left[\sup_{\|\beta - \widetilde{\beta}\| < \delta}\|m(Z_i, \beta) - m(Z_i, \widetilde{\beta})\|^2 \times \sup_h\left|\frac{p(X_i)}{1 - p(X_i)}\right|^2\right]$$

$$+ E_a\left[\sup_{\widetilde{\beta}\in B:\,\|\widetilde{\beta}-\beta_o\|\leq\delta_0} \|m(Z_i,\widetilde{\beta})\|^2\right.$$

$$\left.\times \sup_{\|p(\cdot)-\widetilde{p}(\cdot)\|_{\infty,\omega}<\delta}\left|\frac{p(X_i)}{1-p(X_i)} - \frac{\widetilde{p}(X_i)}{1-\widetilde{p}(X_i)}\right|^2\right]$$

$$\leq E_a\left[\sup_{\|\beta-\widetilde{\beta}\|<\delta}\|m(Z_i,\beta) - m(Y_i,X_i,\widetilde{\beta})\|^2\right] \times \left(\frac{\overline{p}}{1-\overline{p}}\right)^2$$

$$+ E_a\left[\sup_{\widetilde{\beta}\in B\|\widetilde{\beta}-\beta_o\|\leq\delta_0} \|m(Z_i,\widetilde{\beta})\|^2(1+|X_i|^2)^\omega\right]$$

$$\times \sup_{\|p(\cdot)-\widetilde{p}(\cdot)\|_{\infty,\omega}<\delta}\sup_{x\in\mathcal{X}}|[p(x)-\widetilde{p}(x)](1+|x|^2)^{-\omega/2}|^2 \frac{1}{(1-\overline{p})^2}$$

$$\leq const.\delta^{2\epsilon} + const.\delta^2 \qquad \text{for some } \epsilon \in (0,1],$$

where the last inequality is due to our Assumptions 4.2, 4.3 and Proposition B.1(i). In the following we let $N(\varepsilon,\Lambda_c^\gamma(\mathcal{X}),\|\cdot\|_{\infty,\omega})$ denote the $\|\cdot\|_{\infty,\omega}$-covering number of $\Lambda_c^\gamma(\mathcal{X})$ [i.e., the minimal number of $N$ for which there exist $\varepsilon$-balls $\{h:\|h-u_j\|_{\infty,\omega}\leq\varepsilon\}$, $j=1,\ldots,N$ to cover $\Lambda_c^\gamma(\mathcal{X})$]. Then our Assumption 4.1 implies

$$\log N(\delta,\Lambda_c^\gamma(\mathcal{X}),\|\cdot\|_{\infty,\omega}) \leq const.\left(\frac{c}{\delta}\right)^{d_x/\gamma},$$

$$\int_0^1 \sqrt{\log N(\delta,\Lambda_c^\gamma(\mathcal{X}),\|\cdot\|_{\infty,\omega})}\,d\delta < \infty.$$

Thus by applying CLK's Theorem 3, CLK's condition (2.5) is satisfied.

It remains to verify CLK's condition (2.6). First we note

$$\sqrt{n_a}M_n(\beta_o,p_o) = \frac{1}{\sqrt{n_a}}\sum_{i=1}^{n_a} m(Z_i,\beta_o)\frac{p_o(X_i)}{1-p_o(X_i)}$$

$$= \frac{1}{\sqrt{n_a}}\sum_{i=1}^{n}(1-D_i)m(Z_i,\beta_o)\frac{p_o(X_i)}{1-p_o(X_i)}$$

$$= \sqrt{\frac{n}{n_a}} \times \frac{1}{\sqrt{n}}\sum_{i=1}^{n}(1-D_i)m(Z_i,\beta_o)\frac{p_o(X_i)}{1-p_o(X_i)}.$$

Next we notice

$$\sqrt{n_a}\Gamma_2(\beta_o,p_o)[p(\cdot)-p_o(\cdot)]$$

$$= \sqrt{n_a}E_a\left\{\mathcal{E}(X,\beta_o)\frac{p(X)-p_o(X)}{(1-p_o(X))^2}\right\}$$

$$= \sqrt{\frac{n_a}{n}} \frac{1}{1-p} \times \sqrt{n} E\left\{\mathcal{E}(X,\beta_o)\frac{p(X)-p_o(X)}{1-p_o(X)}\right\}.$$

By Lemma B.2 and $n_a/n = 1 - p + o_p(1)$, we obtain

$$\sqrt{n_a}\{M_n(\beta_o, p_o) + \Gamma_2(\beta_o, p_o)[\widehat{p}(\cdot) - p_o(\cdot)]\}$$

$$= \sqrt{\frac{1}{1-p}}$$

$$\times \frac{1}{\sqrt{n}}\sum_{i=1}^n\left\{(1-D_i)m(Z_i,\beta_o)\frac{p_o(X_i)}{1-p_o(X_i)} + \frac{D_i - p_o(X_i)}{1-p_o(X_i)}\mathcal{E}(X,\beta_o)\right\}$$

$$+ o_p(1),$$

thus CLK's condition (2.6) is satisfied. Moreover from the proof of CLK's Theorem 2 we obtain

$$\sqrt{n_a}(\widehat{\beta} - \beta_o) = -(\Gamma_1'W\Gamma_1)^{-1}\Gamma_1'W\sqrt{n_a}\{M_n(\beta_o, p_o)$$

$$+ \Gamma_2(\beta_o, p_o)[\widehat{p}(\cdot) - p_o(\cdot)]\} + o_p(1)$$

$$= -\frac{1-p}{p}(J_\beta^{1\prime}WJ_\beta^1)^{-1}J_\beta^{1\prime}W\sqrt{n_a}\{M_n(\beta_o, p_o)$$

$$+ \Gamma_2(\beta_o, p_o)[\widehat{p}(\cdot) - p_o(\cdot)]\} + o_p(1).$$

Since $\frac{n}{n_a} = \frac{1}{1-p} + o_p(1)$,

$$\sqrt{n}(\widehat{\beta} - \beta_o) = -\frac{1-p}{p}(J_\beta^{1\prime}WJ_\beta^1)^{-1}J_\beta^{1\prime}W\sqrt{n}\{M_n(\beta_o, p_o)$$

$$+ \Gamma_2(\beta_o, p_o)[\widehat{p}(\cdot) - p_o(\cdot)]\} + o_p(1)$$

$$= -(J_\beta^{1\prime}WJ_\beta^1)^{-1}J_\beta^{1\prime}W\frac{1}{p}\frac{1}{\sqrt{n}}$$

$$\times \sum_{i=1}^n\left\{m(Z_i,\beta_o)\frac{[1-D_i]p_o(X_i)}{1-p_o(X_i)} + \frac{D_i - p_o(X_i)}{1-p_o(X_i)}\mathcal{E}(X,\beta_o)\right\}$$

$$+ o_p(1),$$

thus we obtain Theorem 7 after we establish condition (*).

We now apply Assumption 4.6a or 4.6b or 4.6c to verify condition (*). Since

$$M(\beta, p(\cdot)) - M(\beta, p_o(\cdot)) - \Gamma_2(\beta, p_o)[p(\cdot) - p_o(\cdot)]$$

$$= E_a\left\{m(Z,\beta)\left[\frac{p(X)}{1-p(X)} - \frac{p_o(X)}{1-p_o(X)} - \frac{p(X)-p_o(X)}{(1-p_o(X))^2}\right]\right\}$$

$$= E_a\left\{\frac{m(Z,\beta)[p(X)-p_o(X)]}{1-p_o(X)}\left[\frac{1}{1-p(X)} - \frac{1}{1-p_o(X)}\right]\right\}$$

$$= E_a\left\{\frac{\mathcal{E}(X,\beta)[p(X)-p_o(X)]^2}{(1-p(X))(1-p_o(X))^2}\right\},$$

we have under Assumption 3.1,

$$\sup_{\beta\in B:\,\|\beta-\beta_0\|\leq\delta_0}\|M(\beta,p(\cdot))-M(\beta,p_o(\cdot))-\Gamma_2(\beta,p_o)[p(\cdot)-p_o(\cdot)]\|$$

$$=\sup_{\beta\in B:\,\|\beta-\beta_0\|\leq\delta_0}\left\|E_a\left\{\frac{\mathcal{E}(X,\beta)[p(X)-p_o(X)]^2}{(1-p(X))(1-p_o(X))^2}\right\}\right\|$$

$$\leq\frac{1}{(1-\overline{p})^3}E_a\left\{\sup_{\beta\in B:\,\|\beta-\beta_0\|\leq\delta_0}\|\mathcal{E}(X,\beta)\|\times[p(X)-p_o(X)]^2\right\}.$$

If Assumption 4.6a holds, then

$$E_a\left\{\sup_{\beta\in B:\|\beta-\beta_0\|\leq\delta_0}\|\mathcal{E}(X,\beta)\|\times[p(X)-p_o(X)]^2\right\}$$

$$\leq\sup_{\beta\in B:\,\|\beta-\beta_0\|\leq\delta_0}\sup_x\|\mathcal{E}(x,\beta)\|\times E_a\{[p(X)-p_o(X)]^2\}$$

$$\leq const.[\|p(\cdot)-p_o(\cdot)\|_{2,a}]^2.$$

Now Proposition B.1(ii), $k_n=O(n^{d_x/(2\gamma+d_x)})$ and $\gamma>d_x/2$ imply $[\|\widehat{p}(\cdot)-p_o(\cdot)\|_{2,a}]^2=o_p(n^{-1/2})$, hence condition (*) is satisfied.

If Assumption 4.6b holds, then

$$E_a\left\{\sup_{\beta\in B:\,\|\beta-\beta_0\|\leq\delta_0}\|\mathcal{E}(X,\beta)\|\times[p(X)-p_o(X)]^2\right\}$$

$$\leq\left(E_a\left[\sup_{\beta\in B:\,\|\beta-\beta_0\|\leq\delta_0}\|\mathcal{E}(X,\beta)\|^4\right]\right)^{1/4}(E_a\{[p(X)-p_o(X)]^4\})^{1/4}$$

$$\times\sqrt{E_a\{[p(X)-p_o(X)]^2\}}$$

$$\leq const.\times[\|p(\cdot)-p_o(\cdot)\|_{2,a}]^{2-d_x/(4\gamma)}$$

$$\text{for all }\|p(\cdot)-p_o(\cdot)\|_{2,a}=o(1),$$

where the last inequality is due to the following inequalities for any $s\in[\frac{d_x}{4},\gamma)$:

$$(E_a\{[p(X)-p_o(X)]^4\})^{1/4}\leq const.(\|p(\cdot)-p_o(\cdot)\|_{2,a}+\|\nabla^s\{p(\cdot)-p_o(\cdot)\}\|_{2,a}),$$

$$\|\nabla^s\{p(\cdot)-p_o(\cdot)\}\|_{2,a}\leq const.[\|p(\cdot)-p_o(\cdot)\|_{2,a}]^{1-s/\gamma}.$$

Now Proposition B.1(ii), $k_n=O(n^{d_x/(2\gamma+d_x)})$ and $\gamma>3d_x/4$ imply $[\|\widehat{p}(\cdot)-p_o(\cdot)\|_{2,a}]^{2-d_x/(4\gamma)}=o_p(n^{-1/2})$, hence condition (*) is satisfied.

If Assumption 4.6c holds, then

$$E_a\left\{\sup_{\beta\in B:\|\beta-\beta_0\|\leq\delta_0}\|\mathcal{E}(X,\beta)\|\times[p(X)-p_o(X)]^2\right\}$$
$$\leq\sqrt{E_a\left[\sup_{\beta\in B:\|\beta-\beta_0\|\leq\delta_0}\|\mathcal{E}(X,\beta)\|^2\right]}\times\sqrt{E_a\{[p(X)-p_o(X)]^4\}}$$
$$\leq const.\times[\|p(\cdot)-p_o(\cdot)\|_{2,a}]^{2(1-d_x/(4\gamma))}$$

for all $\|p(\cdot)-p_o(\cdot)\|_{2,a}=o(1)$.

Now Proposition B.1(ii), $k_n = O(n^{d_x/(2\gamma+d_x)})$ and $\gamma > d_x$ imply $[\|\widehat{p}(\cdot)-p_o(\cdot)\|_{2,a}]^{2(1-d_x/(4\gamma))}=o_p(n^{-1/2})$, hence condition (*) is satisfied. $\square$

**Acknowledgments.** We thank the co-editor, the associate editor and three referees for very useful suggestions that greatly improved the presentation of the paper. We also thank J. Hahn, J. Ham, G. Imbens, O. Linton, W. Newey, B. Salanié and seminar participants at various universities for comments.

X. Chen  
Department of Economics  
Yale University  
New Haven, Connecticut 06520-8281  
USA  
E-mail: xiaohong.chen@yale.edu

H. Hong  
Department of Economics  
Stanford University  
Stanford, California 94305  
USA  
E-mail: doubleh@stanford.edu

A. Tarozzi  
Department of Economics  
Duke University  
Durham, North Carolina 27708  
USA  
E-mail: taroz@econ.duke.edu